\pdfoutput=1

\documentclass[12pt]{amsart}
\usepackage[left=1in, right=1in, top=1in, bottom =1in]{geometry}

 \usepackage[foot]{amsaddr} 
\usepackage{amsfonts, amsmath, amssymb, amsthm} 
\usepackage[ruled,vlined]{algorithm2e}
\usepackage{float} 
\usepackage{orcidlink} 
\usepackage[export]{adjustbox}
\usepackage{array}
\newcommand*{\affaddr}[1]{#1} 
\newcommand*{\affmark}[1][*]{\textsuperscript{#1}}

\newtheorem{theorem}{Theorem}
\newtheorem{lemma}{Lemma}
\graphicspath{{figures/}}
\usepackage{hyperref}
\hypersetup{breaklinks=true,colorlinks,allcolors=black}
\begin{document}

\title{Regularized Step Directions in Nonlinear Conjugate Gradient Methods}

\author{Cassidy K. Buhler \orcidlink{0000-0003-4157-4273}\affmark[1]}
\author{Hande Y. Benson \orcidlink{0000-0002-5554-9928} \affmark[1]}
\author{David F. Shanno\affmark[2]}

\address{
\affaddr{\affmark[1]Drexel University, Department of Decision Sciences and MIS, Philadelphia, PA, USA}}
\address{\affaddr{\affmark[2]Rutgers University, RUTCOR (Emeritus), New Brunswick, NJ, USA}
}

\email{cassidy.buhler@dragons.drexel.edu}
\email{benson@drexel.edu, shannod@comcast.net}


\begin{abstract}

Conjugate gradient minimization methods (CGM) and their accelerated variants are widely used. We focus on the use of cubic regularization to improve the CGM direction independent of the step length computation. In this paper, we propose the Hybrid Cubic Regularization of CGM, where regularized steps are used selectively. Using Shanno's reformulation of CGM as a memoryless BFGS method, we derive new formulas for the regularized step direction. We show that the regularized step direction uses the same order of computational burden per iteration as its non-regularized version. Moreover, the Hybrid Cubic Regularization of CGM exhibits global convergence with fewer assumptions. In numerical experiments, the new step directions are shown to require fewer iteration counts, improve runtime, and reduce the need to reset the step direction. Overall, the Hybrid Cubic Regularization of CGM exhibits the same memoryless and matrix-free properties, while outperforming CGM as a memoryless BFGS method in iterations and runtime.
\end{abstract}

\maketitle 
\newpage

\section*{Acknowledgements}
Sadly, David F. Shanno passed away in July 2019.  The research documented in this paper was started by Benson and Shanno in 2015 and was presented at the SIAM Optimization Meeting in 2017. Cassidy Buhler joined the research group after Shanno's passing.  This work represents the last project that Benson and Shanno completed in their 20-years of collaboration, and the authors hope that it is a tribute to his legacy and a solid foundation for the next generation of researchers.

The authors would like to thank Drs. M\"{u}ge \c{C}apan, Vasilis Gkatzelis, Chelsey Hill, and Matthew Schneider for their feedback on an earlier version of the paper.  We are especially grateful to Dr. Gkatzelis for conversations on the theoretical results in the paper.

\section{Introduction}
\label{intro}
The unconstrained nonlinear programming problem (NLP) has the form 
\begin{equation} \label{UnconsNLP}
	\min_x f(x)
\end{equation}
where $x \in \mathbb{R}^n$ and $f: \mathbb{R}^n \rightarrow \mathbb{R}$.  We assume that $f:\mathbb{R}^n \rightarrow \mathbb{R}$ is smooth and its Hessian is Lipschitz continuous on at least the set $\{x \in \mathbb{R}^n: f(x) \le f(x_0)\}$.  There are a number of different methods for solving \eqref{UnconsNLP} that are usually classified by the amount of derivative information used.  For instances where function, gradient, or Hessian evaluations may be unavailable or expensive to evaluate, store, or manipulate, the choice of algorithm may be dictated by what is computationally tractable.  Whenever all quantities are readily available, the comparative performance is generally one where there is a trade-off between the number of iterations and the amount of work done per iteration.

In this paper, we focus on conjugate gradient minimization methods (CGM), which are first-order methods for solving \eqref{UnconsNLP}.  These methods are also referred to as {\em nonlinear conjugate gradient methods} in literature to differentiate them from the classical conjugate gradient methods that were outlined in \cite{hestenes1952methods} to solve a linear system of equations. 

It is designed to improve on using only the gradient direction by adding a momentum term: while solving \eqref{UnconsNLP}, CGM generates a sequence of iterates $\{x_k\}$ such that 
\begin{align} \label{seq}
	x_{k+1} &= x_k + \alpha_k \Delta x_k \\
	\Delta x_{k+1} &= -\nabla f(x_{k+1}) + \beta_k \Delta x_k
\end{align}
where $\alpha_k$ is the step length, $\Delta x_k$ is the step direction, and $\beta_k$ is a scalar defined as 

\begin{equation} \label{genBeta}
	\beta_k = \dfrac{ \nabla f(x_{k+1})^T y_k }{y_k^T \Delta x_k}.
\end{equation}
Here, we use the standard notation $y_k := \nabla f(x_{k+1}) - \nabla f(x_k)$.  

If an exact line search is used, under the smoothness assumption, $\alpha_k$ satisfies the first order condition $\Delta x_k^T \nabla f(x_k + \alpha_k \Delta x_k) = 0$, 
or $\Delta x_k^T \nabla f(x_{k+1}) = 0$.
Under the same assumptions, it can also be shown that $\Delta x_k^T \nabla f(x_k) = -\nabla f(x_k)^T \nabla f(x_k)$.  Therefore, with exact line 
search, \eqref{genBeta} can be written as
\begin{equation} \label{exactBeta}
	\beta_k = \dfrac{ \nabla f(x_{k+1})^T y_k }{\nabla f(x_k)^T \nabla f(x_k)}.
\end{equation}
This form of $\beta_k$ gives the Polak-Ribiere formula \cite{PR69}.  If $f$ is quadratic, then \eqref{exactBeta} further reduces to 
\begin{equation} \label{exactBetaQ}
	\beta_k = \dfrac{ \nabla f(x_{k+1})^T \nabla f(x_{k+1})  }{\nabla f(x_k)^T \nabla f(x_k)},
\end{equation}
which is the Fletcher-Reeves formula \cite{fletcher1964function}.

While the momentum term generally yields an improvement over the steepest descent direction, two concerns still remain: 
\begin{enumerate}
    \item The step directions can fail to satisfy the conjugacy condition, and 
    \item The step length calculation can be a bottleneck for runtime as it requires multiple function evaluations.
\end{enumerate}
The first concern can have a significant impact on the number of iterations to reach the solution of \eqref{UnconsNLP}.  In fact, it is shown in \cite{powell1976some} and \cite{crowder1972linear} that CGM as defined by \eqref{seq} and \eqref{genBeta} exhibits a linear rate of convergence unless restarted every $n$ iterations with the steepest descent direction.  Moreover, in addition to restarting the method every $n$ iterations, \cite{powell1977restart} proposed to use a restart whenever the algorithm moves too far from conjugacy, or more precisely when
\begin{equation} \label{PowellCri}
	| \nabla f(x_{k+1})^T \nabla f(x_k) | \ge 0.2 \| \nabla f(x_{k+1}) \|^2.
\end{equation}
In this paper, a restart that occurs after \eqref{PowellCri} will be called a {\em Powell restart}.
We will show in the numerical results section that nearly all of the problems in our test set require at least one Powell restart and, on average, nearly half of all iterations are Powell restart iterations.  Therefore, in order to truly distinguish CGM from steepest descent and improve the rate of convergence, we need a mechanism to improve the step directions.  A formal definition of ``improvement'' will be provided in the next section.

The second concern is about the step length calculation and impacts the amount of effort required per iteration, and, thus, the runtime of the algorithm.
An exact line search seeks to find a step length $\alpha^*$ which solves
\begin{equation} \label{alphamin}
    \min_{\alpha} \Phi(\alpha) = f(x + \alpha \Delta x)
\end{equation}
for a given point $x$ and step direction $\Delta x$.  For a general nonlinear function $f$, this minimization problem in one variable ($\alpha$) is usually ``solved'' using an iterative approach such as bisection or cubic interpolation \cite{davidon1968variance}, even though these approaches cannot find the exact value of the minimizing $\alpha$ within a finite number of iterations.  For specific forms of $f$, such as a strictly convex quadratic function, a formula can be used to 
directly calculate the minimizing $\alpha$ without the need for an iterative approach.  We should also note that at the solution of the exact line search, 
$\Phi'(\alpha) = \nabla f(x + \alpha \Delta x)^T \Delta x = 0$.

By contrast, an inexact line search only seeks to approximately minimize $\Phi(\alpha)$, requiring lower levels of accuracy when the iterates $x_k$ are away from a 
stationary point of $f$.  In order to maintain theoretical guarantees, most inexact line search techniques rely on guaranteeing {\em sufficient descent}, that is,
$f(x) - f(x + \alpha \Delta x)$ must be sufficiently large according to some criterion.  The Armijo criterion is given by 
\[ f(x + \alpha \Delta x) < f(x) + \epsilon_1 \alpha \nabla f(x)^T \Delta x, \]
and it is accompanied by the curvature condition 
\[ -\nabla f(x + \alpha \Delta x)^T \Delta x \le \epsilon_2 \nabla f(x)^T \Delta x, \]
for constants $0 < \epsilon_1 < \epsilon_2 < 1$.  The two conditions together are called the {\em Wolfe conditions}.  The curvature condition can also be modified as
\[ | \nabla f(x + \alpha \Delta x)^T \Delta x| \le \epsilon_2 |\nabla f(x)^T \Delta x| \]
to give the {\em Strong Wolfe conditions}.  

While the two best-known forms of CGM are given by \eqref{exactBeta} and \eqref{exactBetaQ} above, they do rely on using an exact line search. 
This means that when an explicit formula for directly computing the minimizing $\alpha$ is not available, the function $f$ and/or its gradient must be evaluated multiple times for an iterative line search within each iteration of the CGM algorithm.  Doing so can be costly if $n$ is large or if the function evaluation is time consuming.  For machine learning problems, many gradient-based algorithms choose a fixed step length (also referred to as the learning rate) or use an adaptive approach to setting it.   Adagrad \cite{Adagrad}, Adadelta \cite{Adadelta}, RMSprop \cite{RMSprop}, and Adam \cite{Adam} are examples of adaptive learning rate approaches with good performance, but their use does not currently extend to CGM.

In this paper, we propose a cubic regularization variant of CGM, that combines ideas from \cite{BS14} and \cite{benson2018cubic} to selectively use cubic regularization when solving \eqref{UnconsNLP} and from \cite{shanno1978conjugate} to recast CGM as ``memoryless BFGS'' and apply an inexact line search.  The resulting approach is shown to reduce the need for Powell restarts and to (approximately) optimize the step direction without the typical overload of additional function evaluations.  We show in the numerical results that our implementation improves iteration counts and run time on the CUTEst test set \cite{gould2015cutest} and on randomly generated machine learning problems.  As noted in \cite{Gri81} and \cite{BS14}, the approach has an inherent connection to Levenberg-Marquardt's method \cite{Lev44}, \cite{Mar63} and is, therefore, related to the scaled conjugate gradient method proposed in \cite{moller1993scaled} for training neural networks. 
Scaled CGM is the default training function for MATLAB's pattern recognition neural network function, {\tt patternnet} \cite{matlabDoc}.  

In the next section, we introduce the central question of our research and set the vocabulary and notation for the remainder of the paper.  Given the extensive body of literature our work pulls from, we believe that clarifying our goals and our scope into a single framework is necessary to clearly communicate our proposed approach and put our results into perspective.  We then introduce the cubic regularization of CGM, including a review of the memoryless BFGS formulation proposed by Shanno \cite{shanno1978conjugate}, the explicit formulae for the regularized step direction, and the corresponding method for choosing the step length.  Unlike other applications of cubic regularization, such as \cite{CGT1}, \cite{BS14}, \cite{benson2018cubic}, the regularized step in our framework can be computed without significant overhead beyond that required for \eqref{seq} with \eqref{genBeta}.  Moreover, we show that the burden of optimizing the step length can be shifted to optimizing the regularization parameter, which does not require as many function or gradient evaluations as solving \eqref{alphamin}.  In Section 3, we present some theoretical results on computational complexity per iteration, as well as global convergence.  Numerical results are given in Section 4.

{\textbf{Notation:}  Throughout the paper, $\| \cdot \|$ refers to the Euclidean norm.} 

\section{Cubic Regularization for CGM}
\label{sec:1}
As discussed, CGM is designed to be an ``accelerated'' version of steepest descent.  The momentum term for the step direction is obtained cheaply so as not to increase the computational burden over the gradient direction, and the local convergence rate is improved from linear to superlinear.  However, if the step direction moves away from conjugacy and CGM has to be restarted frequently, CGM will behave more like steepest descent.

The goal of this paper is to propose a regularized method to improve step quality within CGM.  Specifically, we aim for this method to
\begin{itemize}
    \item require fewer iterations in computational experiments than its non-regularized version,
    \item exhibit global convergence with fewer assumptions than its non-regularized version,
    \item require the same order of computational burden per iteration as its non-regularized version, and
    \item demonstrate faster overall runtime in computational experiments than its non-regularized version.
\end{itemize}

In this section, we will present our proposed cubic regularization scheme and its related step length rules.  The integration of cubic regularization for CGM will be closer to the approach taken in \cite{BS14}, wherein Benson and Shanno discuss the equivalence between cubic regularization, Levenberg-Marquardt regularization (this equivalence was originally pointed out by Griewank in \cite{Gri81}), trust-region radius control, and the perturbation of the diagonal of the Hessian matrix for line-search approaches based on Newton's method.  Therefore, in this paper, cubic regularization will generally arise in form of a diagonal perturbation to the approximate Hessian matrix.
(Our previous paper on symmetric rank-1 methods with cubic regularization \cite{benson2018cubic} took the approach of modifying the secant equation, which we are not proposing here but will leave for future work.)  Since the formulation of CGM given by \eqref{seq} and \eqref{genBeta} was matrix-free, we will use Shanno's reformulation of CGM as a memoryless BFGS method \cite{shanno1978conjugate}.  We start with a brief review of the reformulation and then present the proposed new approach.

    \subsection{Memoryless BFGS Formulation of CGM}
    \label{subsec:1}
It was shown in \cite{shanno1978conjugate} that a version of CGM is equivalent to a memoryless BFGS method, and we will use that equivalence here to build our cubic regularization approach.  First, as a reminder and to set notation, the direction is calculated as 
\[	\Delta x_k = -\mathbf{H}_k \nabla f(x_k), \]
where $\mathbf{H}_k \approx \left(\nabla^2 f(x_k) \right)^{-1}$, with $\mathbf{H}_0 = \mathbf{I}$ and $\mathbf{H}_{k+1}$ obtained using the BFGS update formula
\begin{equation} \label{BFGS}
	\mathbf{H}_{k+1} = \mathbf{H}_k - \dfrac{\mathbf{H}_k y_k p_k^T + p_k y_k^T \mathbf{H}_k}{p_k^T y_k} + \left( 1 + \dfrac{y_k^T \mathbf{H}_k y_k}{p_k^T y_k} \right) \dfrac{p_k p_k^T}{p_k^T y_k}. 
\end{equation}
Here, $y_k$ is as before and $p_k = \alpha_k \Delta x_k$.  A memoryless BFGS method would mean that the updates are not accumulated, that is, $\mathbf{H}_k$ is replaced by $\mathbf{I}$ in
the update formula and

\begin{equation} \label{mlessBFGS}
	\mathbf{H}_{k+1} = \mathbf{I} - \dfrac{ y_k p_k^T + p_k y_k^T }{p_k^T y_k} + \left( 1 + \dfrac{y_k^T  y_k}{p_k^T y_k} \right) \dfrac{p_k p_k^T}{p_k^T y_k}.
\end{equation}

To derive the equivalence, Shanno \cite{shanno1978conjugate} notes that Perry \cite{doi:10.1287/opre.26.6.1073} expressed \eqref{seq}-\eqref{exactBeta} in matrix form as
\begin{equation} \label{perrymat}
	\Delta x_{k+1} = -\left( \mathbf{I} - \dfrac{p_k y_k^T}{p_k^T y_k} + \dfrac{p_k p_k^T}{p_k^T y_k} \right) \nabla f(x_{k+1}).
\end{equation}
Shanno \cite{shanno1978conjugate} notes that the matrix in this formulation is not symmetric and adds a further correction:
\[
	\Delta x_{k+1} = -\left( \mathbf{I} - \dfrac{p_k y_k^T}{y_k^T p_k} - \dfrac{y_k p_k^T}{y_k^T p_k} + \dfrac{p_k p_k^T}{p_k^T y_k} \right) \nabla f(x_{k+1}).  
\]
Finally, to ensure that a secant condition is satisfied, the last term in the matrix is re-scaled:
\begin{equation} \label{ShannoPerry}
	\Delta x_{k+1} = -\left[ \mathbf{I} - \dfrac{p_k y_k^T}{y_k^T p_k} - \dfrac{y_k p_k^T}{y_k^T p_k} + \left(1 + \dfrac{y_k^T y_k}{p_k^T y_k} \right) 
						\dfrac{p_k p_k^T}{p_k^T y_k} \right] \nabla f(x_{k+1}).
\end{equation}
The matrix term is exactly the formula \eqref{mlessBFGS} for the memoryless BFGS update.

It is important to note here that, unlike in BFGS, we do not need to store a matrix or a series of updates to calculate $\Delta x_{k+1}$ using \eqref{ShannoPerry}.
Multiplying $\nabla f(x_{k+1})$ through the matrix in \eqref{ShannoPerry} simply requires dot-products, scalar-vector multiplications, and vector addition and subtraction.
As such, each update only requires the storage of 3 vectors of length $n$ and $O(n)$ operations.

Furthermore, for an exact line search, \eqref{ShannoPerry} reduces to the Polak-Ribiere formula \eqref{exactBeta}, thereby ensuring that our proposed cubic regularization
approach remains valid in that case as well.
Finally, one advantage of using \eqref{ShannoPerry} for CGM is that the criterion $p_k^T y_k > 0$ is always satisfied, which ensures that the sequence of step directions it produces remain stable and is required for
most proofs of global convergence, such as the one proposed by \cite{Lenard}.

\subsubsection{Initializations and Restarts}
\label{subsubsec:1}
In the first iteration, CGM can be initialized using the gradient.  However, a two-step process based on the self-scaling proposed in \cite{oren1976optimal} has demonstrated better stability and improved iteration counts \cite{shanno1978matrix}. We will use the same initialization scheme so that
\begin{align} \label{OSinit}
    \mathbf{H}_0 &= \mathbf{I} \nonumber\\
    \mathbf{H}_1 &= \frac{p_0^T y_0}{y_0^T y_0} \left( \mathbf{I} - \frac{p_0 y_0^T + y_0 p_0^T}{p_0^T y_0} + \frac{y_0^T y_0}{p_0^T y_0} \frac{p_0 p_0^T}{p_0^T y_0} \right) + \frac{p_0 p_0^T}{p_0^T y_0}
\end{align}

As discussed, CGM is restarted every $n$ iterations (Beale restart) and
when condition \eqref{PowellCri} is satisfied (Powell restart).  The inverse Hessian approximation at the most recent restart iteration $t$ is given by

\begin{equation} \label{restartH}
    \mathbf{H}_t = \displaystyle \dfrac{p_t^T y_t}{y_t^T y_t} \left(
            \mathbf{I} - \dfrac{p_t y_t^T + y_t p_t^T}{p_t^T y_t} 
            + \dfrac{y_t^T y_t}{p_t^T y_t} \dfrac{p_t p_t^T}{p_t^T y_t}
            \right) + \dfrac{p_t p_t^T}{p_t^T y_t},
\end{equation}
which matches the initialization process given by \eqref{OSinit}.

We incorporate these changes into our approach by replacing $\mathbf{H}_k$ in
\eqref{BFGS} with $\mathbf{H}_t$ as given in \eqref{restartH} instead of $\mathbf{I}$.  As such, the formula for the step direction is also modified from \eqref{ShannoPerry} to
\begin{equation} \label{scaledDx}
    \Delta x_{k+1} = -\left[ \mathbf{H}_t - \dfrac{\mathbf{H}_t y_k p_k^T + p_k y_k^T \mathbf{H}_t}{p_k^T y_k} + \left( 1 + \dfrac{y_k^T \mathbf{H}_t y_k}{p_k^T y_k} \right) \dfrac{p_k p_k^T}{p_k^T y_k} \right] \nabla f(x_{k+1}).
\end{equation}
Note that $\mathbf{H}_t \nabla f(x_{k+1})$ and $\mathbf{H}_t y_k$ can be obtained via dot products and scalar-vector multiplications, so there is still no need to compute or store a matrix.  With these modifications, the CGM algorithm can be fully described as Algorithm 1.

\begin{algorithm}
Pick a suitable $x_0$ and $\epsilon > 0$.

Using the two-step initialization in \eqref{OSinit}: 

Let $\Delta x_0 = -\mathbf{H}_0 \nabla f(x_0)$, choose $\alpha$ to approximately solve \eqref{alphamin}, and let $x_1 = x_0 + \alpha \Delta x_0$.  \\

Let $\Delta x_1 = -\mathbf{H}_1 \nabla f(x_1)$, choose $\alpha$ to approximately solve \eqref{alphamin}, and let $x_2 = x_1 + \alpha \Delta x_1$.

Set $k = 2$ and $t = 1$.

\While { $\| \nabla f(x_k) \| > \epsilon$ } {

    \eIf{ $(k-t) \mod n = 0$ (Beale restart) OR \eqref{PowellCri} is satisfied (Powell restart)} {
        $t \leftarrow k$ \\
        $\Delta x_k \leftarrow -\mathbf{H}_t \nabla f(x_k)$, where $\mathbf{H}_t$ is defined by \eqref{restartH}.
	} {
		$\Delta x_k \leftarrow -\mathbf{H}_k \nabla f(x_k)$, using the formula \eqref{scaledDx}.
	}

    Choose $\alpha$ to approximately solve \eqref{alphamin} given $x_k$ and $\Delta x_k$.
    
    $x_{k+1} \leftarrow x_k + \alpha \Delta x_k$.
    
    $k \leftarrow k + 1$.
}
\caption{Memoryless BFGS reformulation of CGM, as given by \cite{shanno1978conjugate}.}
\end{algorithm}

It should be noted that the line search is the same one in \cite{shanno1978conjugate}, which is an inexact line search that requires sufficient decrease of the objective function at each iteration.  The implementation uses cubic interpolation \cite{davidon1968variance} to approximately solve \eqref{alphamin}.

    \subsection{Cubic Regularization for Quasi-Newton Methods}
\label{subsec:2}
The step-direction, $\Delta x$, used by a quasi-Newton method minimizes
\begin{equation} \label{Quad}
	f_N(x_k + \Delta x) := f(x_k) + \nabla f(x_k)^T \Delta x + \dfrac{1}{2} \Delta x^T \mathbf{B}_k \Delta x,
\end{equation}
where $\mathbf{B}_k \approx \nabla^2 f(x_k)$ and $\mathbf{H}_k = \mathbf{B}_k^{-1}$.

To obtain the cubic regularization formula, we also define
\begin{equation} \label{uncons-CubicM}
	f_M(x_k + \Delta x) := f(x_k) + \nabla f(x_k)^T \Delta x + \dfrac{1}{2} \Delta x^T \mathbf{B}_k \Delta x + \dfrac{M}{6} \| \Delta x \|^3,
\end{equation}
where $M$ is the approximation to the Lipschitz constant for $\nabla^2 f(x)$.  The cubic step direction is found by solving the problem
\begin{equation} \label{cube-prob}
\Delta x \in \arg \min_s f_M (x^k + s).
\end{equation}

In \cite{Gri81}, \cite{NP06}, and \cite{CGT1}, it is shown that for sufficiently large $M$, $x_k + \Delta x$ will satisfy an Armijo condition and that a line search is not needed when using cubic regularization within a quasi-Newton method or Newton's method.  In this new framework, we need to control $M$ rather than $\alpha$.  In \cite{CGT1}, the ARC method starts with a sufficiently large value of $M$ that is decreased through the iterations and approaches (or is set to) 0 in a neighborhood of the solution.  In \cite{BS14} and \cite{benson2018cubic}, the authors proposed setting $M=0$ for all iterations where the Hessian or its estimate are positive definite and picking a value of $M$ using iteration-specific data only as needed.   

In order to motivate the selective use of cubic regularization, we observe that determining $\Delta x$ using \eqref{uncons-CubicM}-\eqref{cube-prob} is nontrivial as it involves the solution of an unconstrained NLP.  In \cite{CGT1}, the authors show that it suffices to solve \eqref{cube-prob} only approximately in order to achieve global convergence.  

Moreover, the use of cubic regularization during iterations with negative curvature is based on its equivalence to the Levenberg-Marquardt method. To see the equivalence, let us examine the solution of \eqref{cube-prob}.  Note that the first-order necessary conditions for the optimization problem are
\begin{equation} \label{1st-cubic}
	\nabla f(x_k) + \left( \mathbf{B}_k + \dfrac{M}{2} \| s \| \mathbf{I} \right) s = 0.
\end{equation}
Similarly, the Levenberg-Marquardt method replaces $\mathbf{B}_k$ with $\mathbf{B}_k + \lambda \mathbf{I}$ for a sufficiently large $\lambda > 0$ to satisfy certain descent criteria.  (A Levenberg-Marquardt-based variant of CGM was proposed by \cite{moller1993scaled} as scaled conjugate gradient methods and remains a popular approach for training neural networks.)  Note that the step, $\Delta x$, obtained by solving
\[
	(\mathbf{B}_k + \lambda \mathbf{I}) \Delta x = -\nabla f(x^k),
\]
satisfies \eqref{1st-cubic} when 
\[
	M = \dfrac{2 \lambda}{\| \Delta x \|}.
\]
(Further details of the equivalence are provided in \cite{BS14}.)
Thus, we simply re-interpret the cubic regularization step as coming from a Levenberg-Marquardt regularization of the CGM step with appropriately related values of $\lambda$ and $M$.  (A similar insight was mentioned in \cite{WDE} but not explicitly used.) If this direction is not accepted, then we increase $\lambda$.  

\subsection{Cubic Regularization for CGM}
\label{subsec:3}
We posit here that the theoretical difficulties encountered by CGM can be similarly addressed via the cubic regularization approach, without reducing and potentially even improving
its overall solution time.  In this section, we start by deriving the update formula for an iteration during which cubic regularization is used.  Then, we will show the impact
of using cubic regularization to reduce the need for restarts in the algorithm.  In Section 4, our numerical results show the effectiveness of this approach in practice.

As discussed in the previous section, the use of cubic regularization necessitates that we add a term to the approximate Hessian, that is, compute the step direction $\Delta x_k$ 
using $\mathbf{B}_{k+1} + \lambda \mathbf{I}$, instead of $\mathbf{B}_{k+1}$.  While the regularization is applied to the approximate Hessian itself, the CGM update formula \eqref{BFGS} updates the
inverse of the approximate Hessian, that is $\mathbf{H}_{k+1} = (\mathbf{B}_{k+1})^{-1}$.  As such, in order to compute the regularized step direction, we need to compute 
$(\mathbf{B}_{k+1} + \lambda \mathbf{I})^{-1}$.  In previous papers that use regularization with BFGS or with Newton's method, there is no direct formula for computing this inverse.  As such, 
the solution of multiple linear systems may be required at each iteration until a suitable $\lambda$ value is found, which means that despite improved step directions that reduce the number of iterations, the effort and time per iteration increase, potentially increasing overall solution time as well.

When applying the regularization to the CG update formula, however, we can derive an explicit formula for $(\mathbf{B}_{k+1} + \lambda \mathbf{I})^{-1}$.  This is a significant advantage for CGM, in that the use of cubic regularization does not incur additional computational burden at each iteration.

We start by showing the formulae for $\mathbf{B}_t$ and $\mathbf{B}_{k+1}$.
\begin{equation} \label{Bt}
	\mathbf{B}_t = \left( \hat{H}_t^{-1} \right) = \left( \dfrac{y_t^T y_t}{p_t^T y_t} \right) \left( \mathbf{I} - \dfrac{p_t p_t^T}{p_t^T p_t} + \dfrac{y_t y_t^T}{y_t^T y_t} \right)
\end{equation}
\begin{equation} \label{Bkp1}
	\mathbf{B}_{k+1} = \left( \hat{H}_{k+1}^{-1} \right) = \mathbf{B}_t - \dfrac{\mathbf{B}_t p_k p_k^T \mathbf{B}_t}{p_k^T \mathbf{B}_t p_k} + \dfrac{y_k y_k^T}{p_k^T y_k}
\end{equation}

Next, we will apply the regularization to $\mathbf{B}_{k+1}$ and take its inverse.  To do so, we will first need to compute the inverse of $\mathbf{B}_t + \lambda \mathbf{I}$ (henceforth referred to as $\mathbf{H}_t(\lambda)$):

\begin{align} \label{Htlambda}
    \mathbf{H}_t(\lambda) &= \left( \mathbf{B}_t + \lambda \mathbf{I} \right)^{-1} \nonumber\\
    &= \frac{p_t^T y_t}{c} \mathbf{I} + \frac{ab}{c(\lambda b + a)} p_t p_t^T 
    - \frac{\lambda}{c(\lambda b + a)} y_t y_t^T \nonumber\\
    &\quad - \frac{a}{c(\lambda b + a)} \left( p_t y_t^T + y_t p_t^T \right),
\end{align}

where
\[
	a = \dfrac{y_t^T y_t}{p_t^T p_t}, \quad b = 2\dfrac{y_t^T y_t}{p_t^T y_t} + \lambda, \quad c = y_t^T y_t + \lambda p_t^T y_t.
\]
It is easy to verify that when $\lambda = 0$, \eqref{Htlambda} reduces to \eqref{restartH}.

We can finally write the formula for computing the inverse of $\mathbf{B}_{k+1} + \lambda \mathbf{I}$ (henceforth referred to as $\mathbf{H}_{k+1}$) by repeatedly applying the Sherman-Morrison-Woodbury formula \cite{sherman1949adjustment}:

\begin{align} \label{Hkp1lambda}
    \mathbf{H}_{k+1} (\lambda) &= \left( \mathbf{B}_{k+1} + \lambda \mathbf{I} \right)^{-1} \nonumber\\
    &= \left( \mathbf{B}_t + \lambda \mathbf{I} - \frac{\mathbf{B}_t p_k p_k^T \mathbf{B}_t}{p_k^T \mathbf{B}_t p_k} + \frac{y_k y_k^T}{p_k^T y_k} \right)^{-1} \nonumber\\
    &= \left( \mathbf{B}_t + \lambda \mathbf{I} \right)^{-1} \nonumber\\
    &\quad + \frac{p_k^T y_k + y_k^T \left( \mathbf{B}_t + \lambda \mathbf{I} \right)^{-1} y_k}{d} \left( \mathbf{B}_t + \lambda \mathbf{I} \right)^{-1} (\mathbf{B}_t p_k) (\mathbf{B}_t p_k)^T \left( \mathbf{B}_t + \lambda \mathbf{I} \right)^{-1} \nonumber\\
    &\quad - \frac{p_k^T \mathbf{B}_t p_k - (\mathbf{B}_t p_k)^T \left( \mathbf{B}_t + \lambda \mathbf{I} \right)^{-1} \mathbf{B}_t p_k}{d} \left( \mathbf{B}_t + \lambda \mathbf{I} \right)^{-1} y_k y_k^T \left( \mathbf{B}_t + \lambda \mathbf{I} \right)^{-1} \nonumber\\
    &\quad - \frac{(\mathbf{B}_t p_k)^T \left( \mathbf{B}_t + \lambda \mathbf{I} \right)^{-1} y_k}{d} \left( \mathbf{B}_t + \lambda \mathbf{I} \right)^{-1} (\mathbf{B}_t p_k) y_k^T \left( \mathbf{B}_t + \lambda \mathbf{I} \right)^{-1} \nonumber\\
    &\quad - \frac{y_k^T \left( \mathbf{B}_t + \lambda \mathbf{I} \right)^{-1} \mathbf{B}_t p_k}{d} \left( \mathbf{B}_t + \lambda \mathbf{I} \right)^{-1} y_k (\mathbf{B}_t p_k)^T \left( \mathbf{B}_t + \lambda \mathbf{I} \right)^{-1},
\end{align}
where the denominator $d$ is given by

\begin{align*}
    d =& \left( p_k^T y_k + y_k^T \left( \mathbf{B}_t + \lambda \mathbf{I} \right)^{-1} y_k \right)
    \left( p_k^T \mathbf{B}_t p_k - (\mathbf{B}_t p_k)^T \left( \mathbf{B}_t + \lambda \mathbf{I} \right)^{-1} (\mathbf{B}_t p_k) \right) \\
    &+ \left( (\mathbf{B}_t p_k)^T \left( \mathbf{B}_t + \lambda \mathbf{I} \right)^{-1} y_k \right)^2.
\end{align*}
In order to see that the formula \eqref{Hkp1lambda} consists of a sum of rank-1 updates to $\mathbf{H}_t(\lambda)$, we introduce the following intermediate calculation:
\begin{align*}
    \tilde{p}_k &= (\mathbf{B}_t + \lambda \mathbf{I})^{-1} \mathbf{B}_t p_k \\
    &= \frac{y_t^T y_t}{c} p_k + \frac{\lambda ab p_t^T p_k}{c(\lambda b + a)} p_t
    - \frac{\lambda^2 y_t^T p_k}{c(\lambda b + a)} y_t 
    - \frac{\lambda a y_t^T p_k}{c(\lambda b + a)} p_t - \frac{\lambda a p_t^T p_k}{c(\lambda b + a)} y_t
\end{align*}
and rewrite \eqref{Hkp1lambda} as
\begin{align} \label{Hkp1lambda2}
    \mathbf{H}_{k+1} (\lambda) &= \mathbf{H}_t(\lambda) - \frac{\tilde{p}_k^T y_k}{d} \left( \tilde{p}_k y_k^T \mathbf{H}_t(\lambda) + \mathbf{H}_t(\lambda) y_k \tilde{p}_k^T \right) \nonumber\\
    &\quad + \frac{p_k^T y_k + y_k^T \mathbf{H}_t(\lambda) y_k}{d} \tilde{p}_k \tilde{p}_k^T \nonumber\\
    &\quad - \frac{p_k^T \mathbf{B}_t p_k - p_k^T \mathbf{B}_t \tilde{p}_k}{d} \mathbf{H}_t(\lambda) y_k y_k^T \mathbf{H}_t(\lambda),
\end{align}
and $d$ as
\[
	d = (p_k^T y_k + y_k^T \mathbf{H}_t(\lambda)y_k) (p_k^T \mathbf{B}_t p_k - p_k^T \mathbf{B}_t \tilde{p}_k) + (y_k^T \tilde{p}_k)^2.
\]
When $\lambda = 0$, we have that $\tilde{p}_k = p_k$ and, therefore, \eqref{Hkp1lambda} reduces to \eqref{BFGS}.

\subsection{\texorpdfstring{Setting a Value for {$\lambda$}}{Setting a Value}}
\label{subsec:4}

Now that we know how to calculate a step direction with cubic regularization, we need to answer two questions:
\begin{enumerate}
	\item When do we apply cubic regularization?
	\item When applying cubic regularization, how do we choose a value for $\lambda$?
\end{enumerate}
The answer to the first question determines how we answer the second one. As shown above, the computation of the regularized step direction does not require significantly more effort than the non-regularized one, so it remains to be seen whether being selective with when to apply the regularization is as important for CGM as it was for Newton's method and quasi-Newton methods.  

To start with, we decided to try finding a pair $(\lambda, \alpha)$ which minimizes $f(x_k + \alpha \Delta x)$, where $\Delta x$ was calculated from a regularized step in every iteration.  Our preliminary numerical studies showed that this approach reduced the number of iterations for many of the problems, but it significantly worsened the computational effort per iteration by requiring an iterative approach that simultaneously optimized over two variables.

Instead, using \cite{BS14} as a guide, we pursued the following approach: selectively use cubic regularization whenever the non-regularized step direction failed to satisfy the Powell criterion, was not 
a descent direction, or resulted in a line search failure.  The assumption
here is that cubic regularization ``improves'' the step direction in some sense, so it should be deployed when the step direction needs such improvement.  In our numerical studies, there were few to no instances of failure to obtain a descent direction or line search failure, so we could not reliably assess the impact of cubic regularization.  However, the prevalence of Powell restarts, as will be noted in Section 4, provided a good opportunity to test potential improvements.

When the non-regularized step leads to a Powell restart, we will set $\lambda > 0$ and try a regularized step, with increasing values of $\lambda$ until it no longer results in a Powell restart.  The initial value of $\lambda$ for a Powell restart is computed as
\begin{equation} \label{laminit}
    \lambda = 5 \dfrac{| \nabla f(x_{k+1})^T \nabla f(x_k) |}{\| \nabla f(x_{k+1}) \|^2}.
\end{equation}
and doubled as needed.  For each value of $\lambda$, we will choose a corresponding optimal $\alpha$.  We have also added a safeguard to bound the number of $\lambda$ updates by a constant $U$, which helps the numerical stability and convergence results of the algorithm.  If the number of updates reaches $U$, we perform a restart. However, in our numerical testing, we set $U=5$ and this bound was never invoked.

With all the details complete, we now describe the approach, called {\em Hybrid Cubic Regularization of CGM}, as Algorithm 2. 

\begin{algorithm}
Pick a suitable $x_0$, $U$, and $\epsilon > 0$.

Using the two-step initialization in \eqref{OSinit}: \\

Let $\Delta x_0 = -H_0 \nabla f(x_0)$, choose $\alpha$ to approximately solve \eqref{alphamin}, and let $x_1 = x_0 + \alpha \Delta x_0$.  

Let $\Delta x_1 = -H_1 \nabla f(x_1)$, choose $\alpha$ to approximately solve \eqref{alphamin}, and let $x_2 = x_1 + \alpha \Delta x_1$.

Set $k = 2$ and $t = 1$.

\While { $\| \nabla f(x_k) \| > \epsilon$ } {

    \eIf{ $(k-t) \mod n = 0$ (Beale restart) } {
        $t \leftarrow k$ \\
        $\Delta x_k \leftarrow -\mathbf{H}_t \nabla f(x_k)$, where $\mathbf{H}_t$ is defined by \eqref{restartH}.\\
        Choose $\alpha$ to approximately solve \eqref{alphamin} given $x_k$ and $\Delta x_k$.
	} {
	    \eIf{ \eqref{PowellCri} is satisfied (Cubic Regularization)} {
	        Reset $k \leftarrow k-1$, $u \leftarrow 1$, and initialize $\lambda$ using \eqref{laminit}.
	        
	        $\Delta x_k \leftarrow -\mathbf{H}_k \nabla  f(x_k)$, where $\mathbf{H}_k$ is defined by \eqref{Hkp1lambda2}.
	        
	        Choose $\alpha$ to approximately solve \eqref{alphamin} given $x_k$ and $\Delta x$.
	        
	        \While{\eqref{PowellCri} is satisfied and $u < U$} {
	            $\lambda \leftarrow 2 \lambda$, $u \leftarrow u + 1$.
	            
	            $\Delta x_k \leftarrow -\mathbf{H}_k  \nabla f(x_k)$, where $\mathbf{H}_k$ is defined by \eqref{Hkp1lambda2}.
	        
	            Choose $\alpha$ to approximately solve \eqref{alphamin} given $x_k$ and $\Delta x_k$.
	        }
	        \If{ $u == U$ } {
	            $t \leftarrow k$ \\
	            $\Delta x_k \leftarrow -\mathbf{H}_k  \nabla f(x_k)$, where $\mathbf{H}_k$ is defined by \eqref{restartH}. \\
	            Choose $\alpha$ to approximately solve \eqref{alphamin} given $x_k$ and $\Delta x_k$.
	        }
	    }{
		$\Delta x \leftarrow -\mathbf{H}_k \nabla f(x_k)$, using the formula \eqref{scaledDx}.
		
		Choose $\alpha$ to approximately solve \eqref{alphamin} given $x_k$ and $\Delta x_k$.
		}
	}

    $x_{k+1} \leftarrow x_k + \alpha \Delta x$.
    
    $k \leftarrow k + 1$.
}
\caption{Hybrid Cubic Regularization of CGM, as described in Section 2.3.}
\end{algorithm}
  
\section{Theoretical Results}
\label{sec:2}
As stated at the beginning of Section 2, we had four dimensions to our goal of improving step quality within CGM.  Two of those goals were theoretical in nature:
\begin{itemize}
    \item exhibit global convergence with fewer assumptions than its non-regularized version.
    \item require the same order of computational burden per iteration as its non-regularized version
\end{itemize}
In this section, we will demonstrate that we have achieved both of these goals.  Since we have mentioned the second goal already, we will start with formalizing it first.

\subsection{Computational Burden per Iteration}
\label{subsec:5}
We start by showing that the additional work per iteration performed by Algorithm 2 does not grow with the size of the problem.  
\begin{theorem}
Computational effort per iteration for Algorithm 2 is of the same order as the computational effort per iteration for Algorithm 1.
\end{theorem}
\proof{}
The Beale restart and the non-regularized steps in Algorithm 2 match those of Algorithm 1.  Therefore, we only need to analyze the steps with cubic regularization.  We know that we will need to try at most $U$ different $\lambda$ values to exit the while-loop in the cubic regularization step with a descent direction that satisfies \eqref{PowellCri} or with a restart.  We can also see that all of the components required to compute $\Delta x$ with cubic regularization using \eqref{Hkp1lambda2}, that is, $\mathbf{H}_t(\lambda)$, $\tilde{p}_k$, $\mathbf{B}_t$, and $d$, can all be obtained using vector-vector and scalar-vector operations of the form that does not necessitate the calculation or storage of a matrix to obtain $-\mathbf{H}_{k+1} \nabla f(x_k)$.  The vectors used in these calculations are the same ones as before ($y_t$, $p_t$, $y_k$, and $p_k$), which means that the computational burden per iteration of the new approach remains the same as in Algorithm 1.
\endproof

\subsection{Global Convergence}
\label{subsec:6}

We now show that Algorithm 2 is globally convergent. Our results utilize the convergence results for Algorithm 1 as provided in \cite{Shanno78conv}, the main theorem of which we have included here as Appendix A for completeness.  The assumptions to establish convergence for Algorithm 1 are stated in \cite{Shanno78conv} as follows:
\begin{itemize}
    \item The eigenvalues of the Hessian of $f$ remain uniformly bounded above.
    \item $f(x)$ is bounded below.    
    \item The line search can find a descent direction providing sufficient descent at each step.
\end{itemize}
The first assumption is used to show that the condition number of the Hessian remains uniformly bounded, which leads directly to the proof of the theorem given here in Appendix A.  The last assumption is to ensure that the algorithm does not cycle.  Additionally, there are mild assumptions (on the line search and the initialization/restart scaling) stated throughout \cite{Shanno78conv} to ensure that the Hessian estimate remains positive definite, so we will consider this an assumption as well, even though it is not explicitly stated in the proof included in Appendix A here.

If the Powell restart/cubic regularization step is not invoked, Algorithm 2 is equivalent to Algorithm 1.  Therefore, it suffices to analyze the iterations with cubic regularization.

The key feature of the proof of the convergence result in \cite{Shanno78conv} is that the condition number of $\mathbf{H}_k$ remains bounded above. Since the only change to the algorithm is to replace $\mathbf{H}_k$ with $\mathbf{H}_k(\lambda) = (\mathbf{B}_k + \lambda \mathbf{I})^{-1}$ in the cubic regularization step, we only need to examine the condition number of $\mathbf{H}_k(\lambda)$ in order to use the rest of the proofs given in \cite{Shanno78conv}. 

\begin{lemma}
 Let $\mathbf{W}$ be a symmetric, positive definite matrix,
and let us denote its Euclidean norm condition number by $\kappa(\mathbf{W})$.  Then, for $\lambda > 0$, \[ \kappa( (\mathbf{W}^{-1} + \lambda \mathbf{I})^{-1}) \le \kappa(\mathbf{W}). \]
\end{lemma}
\proof{}
Let $\chi_{max}$ and $\chi_{min}$ be the maximum and minimum eigenvalues of $\mathbf{W}$, respectively.  By the assumptions made in \cite{Shanno78conv}, we know that $\chi_{min} > 0$.  Then, 

\begin{align*}
    \kappa\left( (\mathbf{W}^{-1} + \lambda \mathbf{I})^{-1} \right) = \kappa\left( \mathbf{W}^{-1} + \lambda \mathbf{I} \right) 
    &= \frac{ (1/\chi_{min}) + \lambda }{ (1/\chi_{max}) + \lambda } \\
    &= \frac{ \chi_{max} + \lambda \chi_{min} \chi_{max} }{ \chi_{min} + \lambda \chi_{min} \chi_{max} } \\
    &\le \frac{ \chi_{max} }{ \chi_{min} } \\
    &= \kappa(\mathbf{W})
\end{align*}
\endproof

With Lemma 1 and the comment above, we can invoke Theorem 2 from Appendix A to conclude that Algorithm 2 exhibits global convergence.

While Algorithm 2 only invokes cubic regularization for a Powell restart, it can be modified to include more cases, such as when the search direction fails to be a descent direction, as another reason to invoke it.  Doing so would mean that the assumption on the line search and the assumption of the Hessian estimate remaining positive definite are no longer necessary. Moreover, by Lemma 1, we can potentially relax the upper bound on the condition number of the Hessian.  We will investigate this in future work.  As such, using cubic regularization allows us to relax one or two of the assumptions in the convergence proof of Algorithm 1.

\section{Numerical Results}
\label{sec:3}
We start this section by describing our testing environment and then we will introduce our numerical results on general unconstrained NLPs from the CUTEst test set \cite{gould2015cutest}.

\subsection{Software and hardware}
\label{subsec:7}

The original conjugate gradient method described in Algorithm 1 was implemented by Shanno and Phua in the code \texttt{Conmin} using Fortran4 \cite{shanno1976algorithm}.  We have reimplemented 
conjugate gradient method of \texttt{Conmin} in C and connected it to {\sc AMPL} \cite{FGK93}.  In our C implementation, we have omitted the BFGS method also implemented in the original \texttt{Conmin} distribution, and we will call our new code \texttt{Conmin-CG}.  The software is available for open source download and use \cite{Conmin-CG}. The cubic regularization scheme proposed in this paper as Algorithm 2 was implemented and tested by modifying this software and is also available at the same link. In our numerical testing, we used {\sc ampl} Version 20210226.

\subsection{Parameters}
In our numerical testing, we set threshold criteria $\epsilon = 1 \times 10^{-6}$ for each algorithm and the bound $U = 5$ for Algorithm 2.

\subsection{Test set}
\label{subsec:8}
The test problems were compiled from the CUTEst test set \cite{gould2015cutest} as implemented in {\sc ampl} \cite{Van97c}.  We chose 230
unconstrained problems, which included all of the unconstrained problems available from \cite{Van97c} except for those where the objective function could not be evaluated
at the provided or default initial solution, where the initial solution was a stationary point (0 iterations), or where the objective function was unbounded below.
There are 38 QPs and 192 NLPs in the set.  We will focus on the two groups of problems separately in the analysis below.

\subsection{Do we need Powell restarts?}
\label{subsec:9}
Of the 192 NLPs in the test set, 187 of them require at least one Powell restart.  That is 97.4\% - a significant portion.  In fact, of the 5 problems that do not
use Powell restarts, 4 conclude within the first two iterations without ever reaching a check for the Powell restart and 1
has an objective function that is the weighted sum of a quadratic term and a nonlinear term, where the weight and function values of the nonlinear term are
negligible with respect to the quadratic term (and as such, it is effectively a QP).  It is, therefore, safe to conclude that every NLP of interest requires at least one Powell
restart and focus on these 187 problems in our ongoing analysis.

Of the 187 problems, our main code fails on 30 problems (8 exit due to line search failures, 3 exit due to function evaluation errors, and 19 reach the iteration limit of 10,000).  The remaining 157 problems are reported as solved and show at least one Powell restart.  For each of these problems, we calculated the 
percentage of Powell restart iterations as
\[
	\begin{array}{lcl}
	\phi_i & = & \text{Percentage of Powell restart iterations for Problem } i \\[.5em]
	& = & \dfrac{ \text{Number of iterations with Powell restarts for Problem } i}{\text{Total number of iterations for Problem } i}.
	\end{array}
\]  
The average value of $\phi$ is 45.5\%, and its median is 45.7\%.  Note that we do not test for a Powell restart during an iteration that is already slated for a Beale restart,
so the percentage of iterations that satisfy the Powell restart criterion \eqref{PowellCri} is actually slightly higher at around 55\%.

Given the prevalence of Powell restarts, it may be natural to ask how the algorithm performs without them.  When Powell restarts are disabled, the code still 
converges on the same number of problems (27 of the 30 failures remain the same, 3 get resolved, and there are 3 new failures).  We get the same iteration count on 13
problems, the code performs better with Powell restarts on 101 problems (80\% fewer iterations on average), and the code performs better without Powell restarts on 41 
problems (30\% fewer iterations on average).  Therefore, there is a clear and significant advantage to Powell restarts, but the question remains as to whether or not 
we need a full restart every time the Powell restart criterion is satisfied.

\subsection{The Impact of Cubic Regularization on the Powell Restart Criterion}
\label{subsec:10}
We will now use an example to illustrate the impact of cubic regularization on the satisfaction of the Powell restart criterion \eqref{PowellCri}.  More specifically,
we will examine how increasing values of $\lambda$ impact the value of the fraction
\begin{equation} \label{PFrac}
	\dfrac{| \nabla f(x_{k+1}(\lambda))^T \nabla f(x_k) |}{\| \nabla f(x_{k+1}(\lambda)) \|^2}.
\end{equation}
The example we have chosen is the problem {\em s206} \cite{gould2015cutest}:
\[
	\min_{x_1,x_2} (x_2 - x_1^2)^2 + 100(1-x_1)^2.
\]
Without cubic regularization, \texttt{Conmin} encounters a Powell restart in Iteration 3.  The value of the Powell fraction \eqref{PFrac} is 23.18, which is much larger than
the threshold of 0.2.  If we start to apply the cubic regularization with increasing values of $\lambda$ for this iteration, we get the results shown in Table \ref{tab:powellRestarts}.

\begin{figure}[ht]
\centering
\includegraphics[width = 105mm,valign=m]{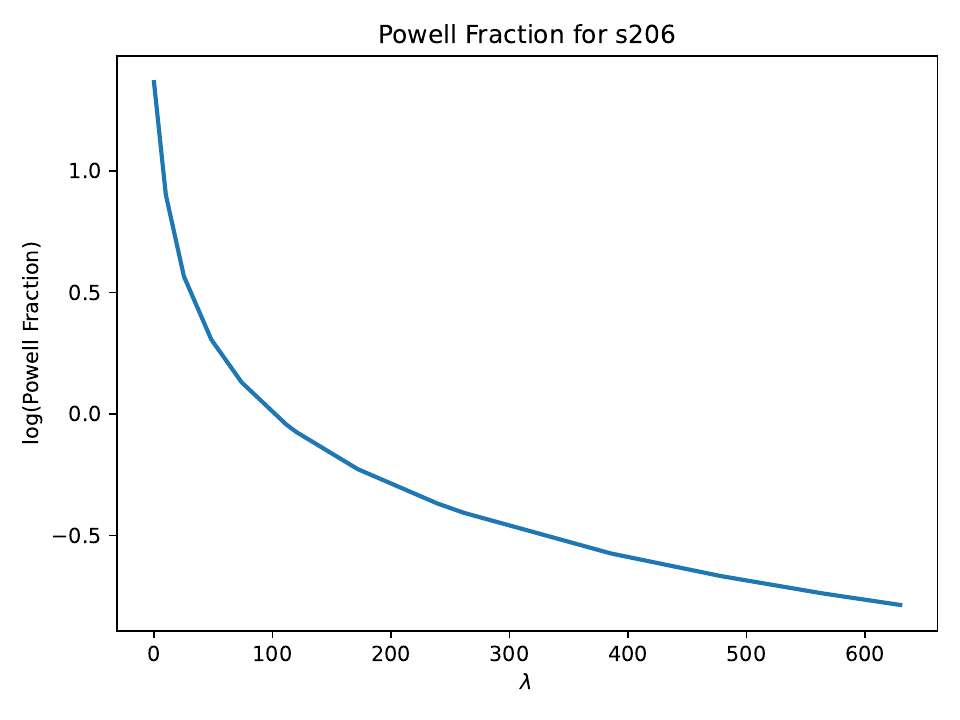}
\caption{$\lambda$ vs. the log of the Powell fraction in Iteration 3 of solving the problem {\em s206}.  The first two iterations of the problem were solved without cubic
regularization. This graph corresponds with Table \ref{tab:powellRestarts}, with more values of $\lambda$ between 0 and 600 to obtain a smooth graph.}
\label{fig:powellRestarts}
\end{figure}

Since the algorithm quickly increases the value of $\lambda$, it is hard to visualize a pattern of decrease with the values given in Table \ref{tab:powellRestarts}. As such, we have
also evaluated the Powell fraction for more values of $\lambda$ between 0 and 600 so that a smooth graph can be produced for Figure \ref{fig:powellRestarts}.

\begin{table}[ht]
    \centering
    \begin{tabular}{|r|r|}
    \hline
        $\lambda$ & Powell Fraction \eqref{PFrac} \\
        \hline
        0.00 & 23.18\\
        115.90 & 0.87\\
        120.44 & 0.84\\
        239.17 & 0.43\\
        478.35 & 0.22\\
        629.78 & 0.16\\
        \hline
    \end{tabular}
\caption{$\lambda$ vs. the log of the Powell fraction in Iteration 3 of solving the problem {\em s206}.  The first two iterations of the problem were solved without cubic
regularization.}
\label{tab:powellRestarts}
\end{table}

\begin{figure}[ht]
    \centering
    \includegraphics[width = 80mm]{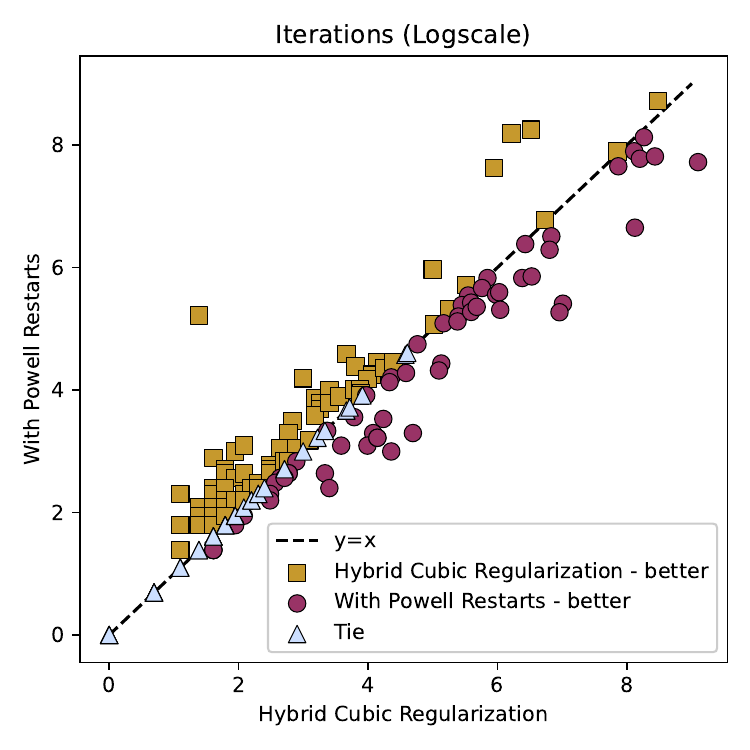} 
    \includegraphics[width = 80mm]{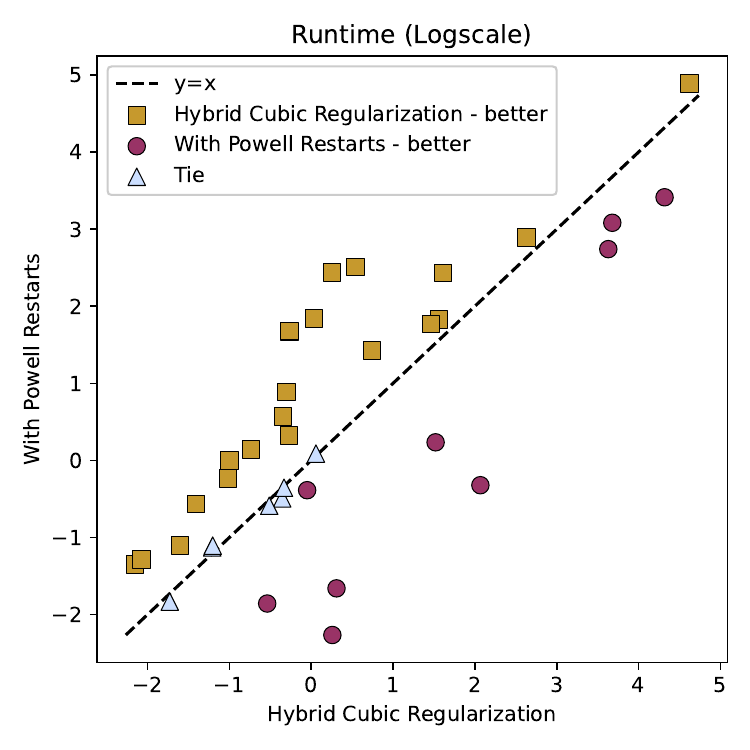}
        \caption[\texttt{Conmin-CG} results on CUTEst test set]{Pairwise comparisons of iterations and runtimes for \texttt{Conmin-CG} with Powell restarts and with hybrid cubic regularization.  The iterations comparison was conducted on 180 out of the 230 unconstrained problems we solved from the CUTEst test set, and the runtimes comparison was conducted on 36 problems on which both solvers exhibited runtimes of at least 0.1 CPU seconds. Yellow squares denote the problems where the code with hybrid cubic regularization outperforms the code with Powell restarts, purple dots represent the opposite relationship, and blue triangles represent a tie (or runtimes within 0.1 of each other). The dotted black line is $y=x$.}
    \label{fig:cutestScatterplot}
\end{figure}

\begin{figure}[ht]
    \centering
    \includegraphics[width = 105mm]{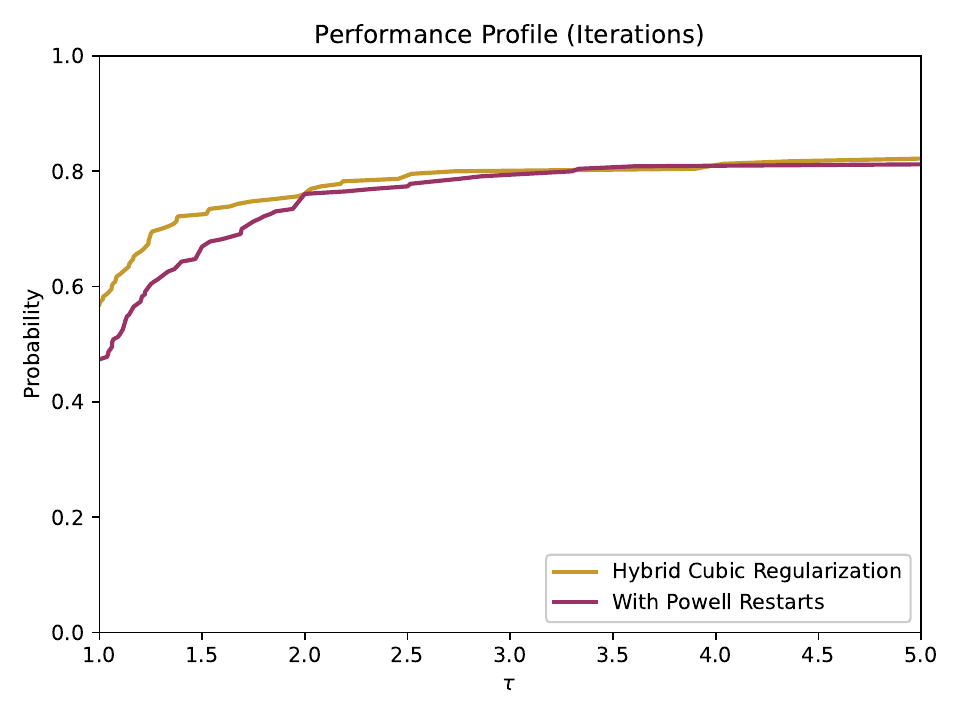}
    \includegraphics[width = 105mm]{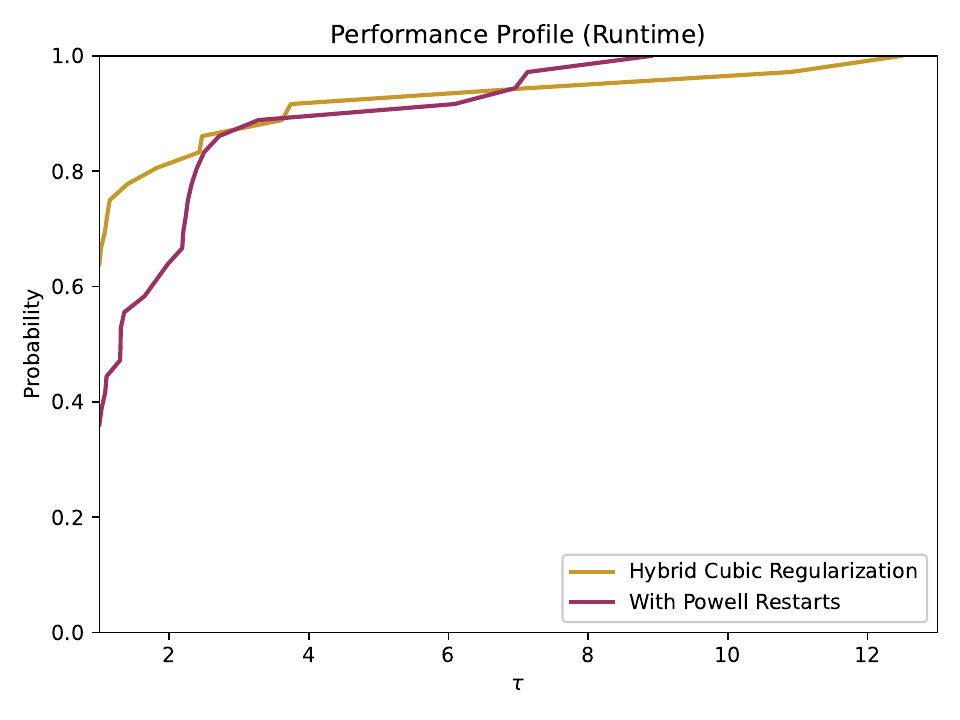}
    \caption{Performance profiles of the iterations and runtime results from CUTEst test set. The iterations comparison was conducted on all 230 unconstrained problems from the CUTEst test set, and the runtimes comparison was conducted on 36 problems on which both solvers exhibited runtimes of at least 0.1 CPU seconds. }
    \label{fig:cutestPerformanceProfiles}
\end{figure}

\subsection{Numerical comparison}
\label{subsec:11}
We start our numerical comparison with general problems of the form \eqref{UnconsNLP} from the CUTEst test set.  Detailed results on these problems are provided in Tables \ref{tab:cutestResults1}-\ref{tab:cutestResults6} in the Appendix.  The detailed results include the number of iterations, runtime in CPU seconds, and the objective function value at the reported solution for each algorithm.  The first, labeled ``With Powell Restarts,'' is the algorithm implemented in \texttt{Conmin-CG} and previously presented as Algorithm 1.  The second, labeled ``No Powell Restarts,'' is a modified version of Algorithm 1 with the check for Powell restarts removed, so that only Beale restarts are performed.  This algorithm is included in the tables to support the results provided in Section 4.3.  The last algorithm, labeled ``Hybrid Cubic,'' implements Algorithm 2, which uses a hybrid approach where cubic regularization is only invoked when the Powell restart criterion \eqref{PowellCri} is satisfied.

The results in Tables \ref{tab:cutestResults1}-\ref{tab:cutestResults6} and Figures \ref{fig:cutestScatterplot} and \ref{fig:cutestPerformanceProfiles} show that hybrid cubic regularization improves the number of iterations and the runtime on the CUTEst test set.  
\begin{itemize}
    \item For overall success, we have that ``With Powell Restarts'' and ``Hybrid Cubic'' each solve 190 of the 230 problems, a rate of 82.6\%.  This includes 180 jointly solved problems, 10 problems solved by ``With Powell Restarts'' only, and 10 problems solved by ``Hybrid Cubic'' only.
    
    \item For iterations, we see in \ref{fig:cutestScatterplot} that ``Hybrid Cubic'' outperforms ``With Powell Restarts'' on the jointly solved problems.  The graph on the left shows a scatterplot where each point represents one problem, with the coordinates equaling iteration counts by the two codes.  (The line $y=x$ is included to help our assessment.)  In this graph, there are 81 yellow squares representing instances where ``Hybrid Cubic'' had fewer iterations, 59 purple dots where ``With Powell Restarts'' had fewer iterations, and 40 blue triangles where they had the same number of iterations.  That means on 121 of the 180 jointly solved problems (67.2\%), ``Hybrid Cubic'' exhibits the same or fewer number of iterations as ``With Powell Restarts.''  The performance profile for iterations, as shown in \ref{fig:cutestPerformanceProfiles}, indicates that ``Hybrid Cubic'' outperforms ``With Powell Restarts'' on all 230 problems of the test set.

    \item It is interesting to note that the improvement becomes slightly more pronounced if we use a typical iteration limit of 1,000 (instead of 10,000).  In that case, there are 166 jointly solved problems.  ``Hybrid Cubic'' performs fewer iterations on 76, ``With Powell Restarts'' performs fewer iterations on 50, and the two codes perform the same number of iterations on 40.  That means on 116 of the 166 jointly solved problems, or 69.9\%, ``Hybrid Cubic'' exhibits the same or fewer number of iterations as ``With Powell Restarts.''
    \item There were 14 jointly solved problems where at least one solver performed over 1,000 iterations.  On these instances, there is no clearly observed pattern, as such a high iteration typically indicates that the problem is severely ill-conditioned or that $f(x)$ is highly nonconvex.  One example is the problem {\em watson}, with $n=31$, and where 
    \[
        f(x) = \sum_{i=1}^{29} \left( 
                    \sum_{j=2}^{n} (j-1)\left(\dfrac{i}{29}\right)^{j-2}x_j - 
                    \left( \sum_{j=1}^{n} \left(\dfrac{i}{29}\right)^{j-2}x_j \right)^2 - 1
                \right)^2 + x_1^2 + (x_2 - x_1^2 - 1)^2.
    \]
    On this problem, ``With Powell Restarts'' performs 2248 and ``Hybrid Cubic'' performs 8919 iterations to reach a solution.  Similar behavior is observed on {\em s371}, which is identical to {\em watson} but with $n=9$.  However, on {\em dixmaani}-{\em dixmaanl}, where $m = 1,000$, $n=3,000$, and 
    \[
        f(x) = 1.0 
            + \sum_{i=1}^{n} \hat{a} \left(\dfrac{i}{n}\right)^2 x_i^2 
            + \sum_{i=1}^{n-1} \hat{b} x_i^2 (x_{i+1} + x_{i+1}^2)^2
            + \sum_{i=1}^{2m} \hat{c} x_i^2 x_{i+m}^4
            + \sum_{i=1}^{m} \hat{d} \left(\dfrac{i}{n}\right)^2 x_i x_{i+2m}
    \]
    for different values of $\hat{a}$, $\hat{b}$, $\hat{c}$, and $\hat{d}$.
    On these problems, ``With Powell Restarts'' performs 2,000-6,000 iterations, which is 1,000-3,000 more than ``Hybrid Cubic.''
    \item For the runtimes comparisons, we have taken a subset of the jointly solved problems, namely those that were solved in 0.1 or more CPU seconds by both solvers.  (This is to ensure a fair comparison, as smaller runtimes can be easily influenced by other processes running on the same machine and/or exhibit very small differences.) The resulting set consisted of 36 problems, of which ``With Powell Restarts'' was faster on 13 and ``Hybrid Cubic'' was faster on 23.  This means that ``Hybrid Cubic'' resulted in faster runtimes on 63.4\% of the problems with runtimes of at least 0.1 CPU seconds.  The runtime results shown in Figures \ref{fig:cutestScatterplot} and \ref{fig:cutestPerformanceProfiles} support this conclusion.
    \item For runtime comparison, it may be a good idea to consider small differences as a tie. If runtimes within 0.1 CPU seconds of each other are considered a tie, then ``With Powell Restarts'' was faster on 9, ``Hybrid Cubic'' was faster on 20, and we considered 7 instances as a tie.  
\end{itemize}

\section{Conclusion}
\label{sec:4}
Our goal in this paper was to incorporate cubic regularization selectively into a CGM framework so that the resulting approach would
\begin{itemize}
    \item require fewer iterations in computational experiments than its non-regularized version,
    \item exhibit global convergence with fewer assumptions than its non-regularized version,
    \item require the same order of computational burden per iteration as its non-regularized version, and
    \item demonstrate faster overall runtime in computational experiments than its non-regularized version.
\end{itemize}
Our global convergence results in Section 3 and our numerical results in Section 4 showed that we were able to attain this goal fully.
    
We have some additional tasks to explore in future work.  In our current implementation, we find the optimal $\alpha$. However, we wish to explore how these results may be affected with a fixed step size. In addition, we hope to extend our work to incorporate subgradients so that we can solve nondifferentiable problems. Finally, we plan to apply CGM towards solving machine learning problems. Our framework is related to a common neural network solver, Scaled Conjugate Gradient \cite{moller1993scaled}. Thus, our work in machine learning will include implementing Hybrid Cubic Regularization of CGM as a solver for neural networks.

\bibliographystyle{plainurl}
\bibliography{ref.bib}

\begin{thebibliography}{10}

\bibitem{benson2018cubic}
Hande~Y Benson and David~F Shanno.
\newblock Cubic regularization in symmetric rank-1 quasi-{N}ewton methods.
\newblock {\em Mathematical Programming Computation}, 10(4):457--486, 2018.

\bibitem{BS14}
H.Y. Benson and D.F. Shanno.
\newblock Interior-point methods for nonconvex nonlinear programming: Cubic regularization.
\newblock {\em Computational Optimization and Applications}, 58, 6 2014.

\bibitem{Conmin-CG}
Cassidy~K Buhler.
\newblock {C}onmin-{CG}, 2024.
\newblock URL: \url{https://github.com/cassiebuhler/conminCG}, \href {https://doi.org/10.5281/zenodo.13315592} {\path{doi:10.5281/zenodo.13315592}}.

\bibitem{CGT1}
C.~Cartis, N.I.M Gould, and Ph.L. Toint.
\newblock Adaptive cubic regularisation methods for unconstrained optimization. {P}art {I}: {M}otivation, convergence and numerical results.
\newblock {\em Mathematical Programming}, 127:245--295, 2011.
\newblock \href {https://doi.org/10.1007/s10107-009-0286-5} {\path{doi:10.1007/s10107-009-0286-5}}.

\bibitem{crowder1972linear}
Harlan Crowder and Philip Wolfe.
\newblock Linear convergence of the conjugate gradient method.
\newblock {\em IBM Journal of Research and Development}, 16(4):431--433, 1972.

\bibitem{davidon1968variance}
William~C Davidon.
\newblock Variance algorithm for minimization.
\newblock {\em The Computer Journal}, 10(4):406--410, 1968.

\bibitem{Adagrad}
John~C Duchi, Elad Hazan, and Yoram Singer.
\newblock Adaptive subgradient methods for online learning and stochastic optimization.
\newblock {\em Journal of Machine Learning Research}, 12:2121--2159, 2011.

\bibitem{fletcher1964function}
Reeves Fletcher and Colin~M Reeves.
\newblock Function minimization by conjugate gradients.
\newblock {\em The Computer Journal}, 7(2):149--154, 1964.

\bibitem{FGK93}
R.~Fourer, D.M. Gay, and B.W. Kernighan.
\newblock {\em {AMPL}: A Modeling Language for Mathematical Programming}.
\newblock Scientific Press, 1993.

\bibitem{gould2015cutest}
Nicholas~IM Gould, Dominique Orban, and Philippe~L Toint.
\newblock {CUTEst}: {A} constrained and unconstrained testing environment with safe threads for mathematical optimization.
\newblock {\em Computational Optimization and Applications}, 60:545--557, 2015.

\bibitem{Gri81}
Andreas Griewank.
\newblock The modification of {N}ewton's method for unconstrained optimization by bounding cubic terms.
\newblock Technical Report NA/12, University of Cambridge, 1981.

\bibitem{hestenes1952methods}
Magnus~Rudolph Hestenes and Eduard Stiefel.
\newblock {\em Methods of Conjugate Gradients for Solving Linear Systems}, volume~49.
\newblock NBS, 1952.

\bibitem{RMSprop}
Geoffrey Hinton, Nitish Srivastava, and Kevin Swersky.
\newblock Neural networks for machine learning lecture 6a overview of mini-batch gradient descent., 2012.
\newblock Retrieved April 27, 2021, from \url{https://www.cs.toronto.edu/~hinton/coursera/lecture6/lec6.pdf}.

\bibitem{Adam}
Diederik~P. Kingma and Jimmy Ba.
\newblock Adam: A method for stochastic optimization.
\newblock In {\em Proceedings of the International Conference on Learning Representations ({ICLR})}, 2015.

\bibitem{Lenard}
M.L. Lenard.
\newblock Practical convergence conditions for unconstrained optimization.
\newblock {\em Mathematical Programming}, 4:309–323, 1973.

\bibitem{Lev44}
K.~Levenberg.
\newblock A method for the solution of certain problems in least squares.
\newblock {\em Quarterly of Applied Mathematics}, 2:164--168, 1944.

\bibitem{Mar63}
D.~Marquardt.
\newblock An algorithm for least-squares estimation of nonlinear parameters.
\newblock {\em SIAM Journal on Applied Mathematics}, 11:431--441, 1963.

\bibitem{matlabDoc}
{MathWorks}.
\newblock Choose a multilayer neural network training function.
\newblock MathWorks Documentation.
\newblock Retrieved April 27, 2021.
\newblock URL: \url{https://www.mathworks.com/help/deeplearning/ug/choose-a-multilayer-neural-network-training-function.html}.

\bibitem{moller1993scaled}
Martin~Fodslette M{\o}ller.
\newblock A scaled conjugate gradient algorithm for fast supervised learning.
\newblock {\em Neural Networks}, 6(4):525--533, 1993.

\bibitem{NP06}
Yu. Nesterov and B.T. Polyak.
\newblock Cubic regularization of {N}ewton method and its global performance.
\newblock {\em Mathematical Programming}, 108:177--205, 2006.
\newblock \href {https://doi.org/10.1007/s10107-006-0706-8} {\path{doi:10.1007/s10107-006-0706-8}}.

\bibitem{oren1976optimal}
Shmuel~S Oren and Emilio Spedicato.
\newblock Optimal conditioning of self-scaling variable metric algorithms.
\newblock {\em Mathematical Programming}, 10(1):70--90, 1976.

\bibitem{doi:10.1287/opre.26.6.1073}
Avinoam Perry.
\newblock Technical note—{A} modified conjugate gradient algorithm.
\newblock {\em Operations Research}, 26(6):1073--1078, 1978.
\newblock \href {https://doi.org/10.1287/opre.26.6.1073} {\path{doi:10.1287/opre.26.6.1073}}.

\bibitem{PR69}
E.~Polak and G.~Rib\`{i}ere.
\newblock Note sur la convergence de m{\'e}thodes de directions conjugu{\'e}es.
\newblock {\em Revue fran{\c{c}}aise d'informatique et de recherche op{\'e}rationnelle. S{\'e}rie rouge}, 16:35--43, 1969.

\bibitem{powell1976some}
Michael James~David Powell.
\newblock Some convergence properties of the conjugate gradient method.
\newblock {\em Mathematical Programming}, 11(1):42--49, 1976.

\bibitem{powell1977restart}
Michael James~David Powell.
\newblock Restart procedures for the conjugate gradient method.
\newblock {\em Mathematical Programming}, 12(1):241--254, 1977.

\bibitem{shanno1978conjugate}
David~F Shanno.
\newblock Conjugate gradient methods with inexact searches.
\newblock {\em Mathematics of Operations Research}, 3(3):244--256, 1978.

\bibitem{shanno1976algorithm}
David~F Shanno and Kang~Hoh Phua.
\newblock Algorithm 500: Minimization of unconstrained multivariate functions [e4].
\newblock {\em ACM Transactions on Mathematical Software (TOMS)}, 2(1):87--94, 1976.

\bibitem{shanno1978matrix}
David~F Shanno and Kang~Hoh Phua.
\newblock Matrix conditioning and nonlinear optimization.
\newblock {\em Mathematical Programming}, 14:149--160, 1978.

\bibitem{Shanno78conv}
D.F. Shanno.
\newblock On the convergence of a new conjugate gradient algorithm.
\newblock {\em SIAM Journal on Numerical Analylsis}, 15(6):1247–1257, 1978.

\bibitem{sherman1949adjustment}
Jack Sherman and Winifred~J Morrison.
\newblock Adjustment of an inverse matrix corresponding to changes in the elements of a given column or a given row of the original matrix.
\newblock In {\em Annals of Mathematical Statistics}, volume 20(4), pages 621--621, 1949.

\bibitem{Van97c}
Robert~J Vanderbei.
\newblock {AMPL} models, 1997.
\newblock URL: \url{https://vanderbei.princeton.edu/ampl/nlmodels/cute/index.html}.

\bibitem{WDE}
M.~Weiser, P.~Deuflhard, and B.~Erdmann.
\newblock Affine conjugate adaptive {N}ewton methods for nonlinear elastomechanics.
\newblock {\em Optimization Methods and Software}, 22(3):413--431, 6 2007.

\bibitem{Adadelta}
Matthew~D Zeiler.
\newblock Adadelta: An adaptive learning rate method.
\newblock In {\em Proceedings of the International Conference on Machine Learning (ICML)}, volume~28, pages 105--112, 2012.
\newblock URL: \url{https://www.jmlr.org/proceedings/papers/v28/zeiler13.pdf}.

\end{thebibliography}

\section{Statements and Declarations}

\subsection{Funding}
The authors declare that no funds, grants, or other support were received during the preparation of this manuscript.

\subsection{Competing Interests}
The authors have no relevant financial or non-financial interests to disclose.

\subsection{Author Contributions}

Hande Benson and David Shanno contributed to the study conception, and all authors contributed to the design. Material preparation, data collection and analysis were performed by Cassidy Buhler and Hande Benson. The first draft of the manuscript was written by Hande Benson and David Shanno, and all authors commented on previous versions of the manuscript. Cassidy Buhler and Hande Benson read and approved the final manuscript.  Approval of the final manuscript was also provided by a family representative of David Shanno after his passing.

\subsection{Data Availability}
The models analyzed during the current study are available in AMPL from \url{https://vanderbei.princeton.edu/ampl/nlmodels/cute/index.html}. 
 These models have been converted to AMPL from SIF and originated from the CUTEst repository, \url{https://github.com/ralna/CUTEst}.

\subsection{Code availability}
    The software \texttt{Conmin-CG}\cite{Conmin-CG} is open source and available for download, \url{https://doi.org/10.5281/zenodo.13315592}. 

\newpage
\section{Appendix}
\subsection{\texorpdfstring{Global Convergence Results from \cite{Shanno78conv}}{Global Convergence Results}}

We showed global convergence in Section \ref{subsec:6} using Theorem 7 from \cite{Shanno78conv}. This theorem is given below for completeness.
\begin{theorem}
        Let $f(x)$ satisfy 
    
    \begin{equation}
    u^TG(x)u \leq m ||u||^2 \quad \textrm{and} \quad f(x) \geq L,
    \end{equation} 
where $u$ is an arbitrary vector in $\mathbb{R}^n$, $G(x) = \nabla^2 f(x)$, $0 < m < \infty$, and $L > -\infty$.
Then, for Algorithm 1, if $\alpha_k$ satisfies
\[
\begin{array}{l}
    p_k^T y_k \geq (-p_k^T \nabla f(x_k)) \epsilon_1, \quad  0 < \epsilon_1 <1 \\
    f(x_{k+1}) - f(x_k) \leq \epsilon_2 \nabla f(x_k)^T p_k, \quad 0 < \epsilon_2 <1 \\
\end{array}
\]
at each step, 
then \[\lim_{k \to \infty} ||p_k|| = 0 \implies \lim \inf_{k \to \infty}||\nabla f(x_k)|| = 0.\]
\end{theorem}


\newpage
\subsection{Detailed Numerical Results}

\begin{table}[b]
\centering
\resizebox{\linewidth}{!}{
\begin{tabular}{|l|r|r|r|r|r|r|r|r|r|r| }
 \hline
  \multicolumn{2}{|c|}{} & \multicolumn{3}{|c|}{With Powell Restarts} &\multicolumn{3}{|c|}{No Powell Restarts} &\multicolumn{3}{|c|}{Hybrid Cubic} \\
 \hline
 Name&n&Iter&Time&f(x*)&Iter&Time&f(x*)&Iter&Time&f(x*)\\
 \hline
aircrftb&5&50&$<$0.1&3.1E-16&53&$<$0.1&3.8E-15&53&$<$0.1&5.2E-13\\
allinitu&4&20&$<$0.1&5.7E+00&9&$<$0.1&5.7E+00&7&$<$0.1&5.7E+00\\
arglina*&100&1&$<$0.1&1.0E+02&1&$<$0.1&1.0E+02&1&$<$0.1&1.0E+02\\
arglinb*&10&1&$<$0.1&4.6E+00&1&$<$0.1&4.6E+00&1&$<$0.1&4.6E+00\\
arglinc*&8&1&$<$0.1&6.1E+00&1&$<$0.1&6.1E+00&1&$<$0.1&6.1E+00\\
arwhead&5000&8&2.6E-01&-9.6E-10&9&$<$0.1&-1.4E-09&4&1.2E-01&-2.7E-09\\
bard&3&17&$<$0.1&8.2E-03&16&$<$0.1&8.2E-03&15&$<$0.1&8.2E-03\\
bdexp&5000&6&2.5E-01&8.0E-04&3&$<$0.1&7.3E-121&3&$<$0.1&7.3E-121\\
bdqrtic&1000&(E)&&&(E)&&&438&4.2E+00&4.0E+03\\
beale&2&11&$<$0.1&8.9E-18&11&$<$0.1&5.5E-19&10&$<$0.1&6.2E-21\\
biggs3&3&14&$<$0.1&1.7E-10&16&$<$0.1&4.5E-13&28&$<$0.1&1.1E-11\\
biggs5&5&84&$<$0.1&5.7E-03&90&$<$0.1&5.7E-03&169&$<$0.1&5.7E-03\\
biggs6&6&62&$<$0.1&3.4E-07&110&$<$0.1&7.2E-09&76&$<$0.1&3.7E-06\\
box2&2&5&$<$0.1&4.2E-15&7&$<$0.1&1.7E-14&5&$<$0.1&3.5E-14\\
box3&3&9&$<$0.1&4.4E-12&10&$<$0.1&2.4E-11&12&$<$0.1&2.9E-11\\
bratu1d&1001&(E)&&&(E)&&&(IL)&&\\
brkmcc&2&4&$<$0.1&1.7E-01&5&$<$0.1&1.7E-01&5&$<$0.1&1.7E-01\\
brownal&10&7&$<$0.1&4.6E-16&7&$<$0.1&1.8E-15&5&$<$0.1&4.0E-15\\
brownbs&2&8&$<$0.1&2.1E-14&11&$<$0.1&2.4E-22&6&$<$0.1&2.2E-13\\
brownden&4&21&$<$0.1&8.6E+04&37&$<$0.1&8.6E+04&14&$<$0.1&8.6E+04\\
broydn7d&1000&339&3.3E-01&4.0E+02&332&2.0E-01&4.0E+02&345&2.0E-01&4.0E+02\\
brybnd&5000&12&3.3E-01&1.5E-12&14&3.0E-01&4.1E-12&13&3.0E-01&4.1E-13\\
chainwoo&1000&167&$<$0.1&1.0E+00&380&1.7E-01&4.6E+00&217&$<$0.1&1.0E+00\\
chnrosnb&50&218&$<$0.1&1.0E-13&258&$<$0.1&3.1E-14&232&$<$0.1&1.1E-13\\
cliff&2&(E)&&&(E)&&&(IL)&&\\
clplatea&4970&877&5.9E+00&-1.3E-02&1306&8.4E+00&-1.3E-02&840&4.3E+00&-1.3E-02\\
clplateb&4970&538&6.2E+00&-7.0E+00&680&4.4E+00&-7.0E+00&900&4.8E+00&-7.0E+00\\
clplatec&4970&(IL)&&&(IL)&&&(IL)&&\\
cosine&10000&7&7.9E-01&-1.0E+04&(E)&&&6&3.6E-01&-1.0E+04\\
cragglvy&5000&80&1.3E+00&1.7E+03&(E)&&&45&4.6E+00&1.7E+03\\
cube&2&14&$<$0.1&2.7E-16&17&$<$0.1&1.1E-17&16&$<$0.1&2.1E-20\\
curly10&10000&(IL)&&&(IL)&&&(IL)&&\\
curly20&10000&(IL)&&&(IL)&&&(IL)&&\\
curly30&10000&(E)&&&(E)&&&(IL)&&\\
 deconvu&51&269&$<$0.1&3.7E-10&321&$<$0.1&3.7E-10&413&$<$0.1&7.9E-08\\
denschna&2&7&$<$0.1&1.2E-19&7&$<$0.1&7.1E-15&7&$<$0.1&7.9E-15\\
denschnb&2&5&$<$0.1&1.2E-14&6&$<$0.1&9.7E-17&5&$<$0.1&2.3E-24\\
\hline
  \end{tabular}}
  \caption{Numerical results on the unconstrained problems from the {CUTE}st test set \cite{gould2015cutest}. Problem names that end in an asterisk are quadratic programming problems. $n$ is the number of variables in the problem, {\em Iter} is the iteration count, {\em Time} is the run time in CPU seconds, and $f(x^*)$ is the objective value at the reported solution.  (IL) and (E) denote that the solver reached its iteration limit and exited with an error, respectively.}
  \label{tab:cutestResults1}
\end{table}

\begin{table}[ht]
\centering
\resizebox{\linewidth}{!}{
\begin{tabular}{|l|r|r|r|r|r|r|r|r|r|r| }
 \hline
  \multicolumn{2}{|c|}{} & \multicolumn{3}{|c|}{With Powell Restarts} &\multicolumn{3}{|c|}{No Powell Restarts} &\multicolumn{3}{|c|}{Hybrid Cubic} \\
 \hline
 Name&n&Iter&Time&f(x*)&Iter&Time&f(x*)&Iter&Time&f(x*)\\
 \hline
denschnc&2&10&$<$0.1&2.2E-17&11&$<$0.1&1.2E-16&8&$<$0.1&1.9E-19\\

denschnd&3&16&$<$0.1&3.7E-11&33&$<$0.1&7.8E-10&12&$<$0.1&2.4E-09\\
denschne&3&(E)&&&(E)&&&(IL)&&\\
denschnf&2&6&$<$0.1&3.9E-21&8&$<$0.1&3.0E-16&4&$<$0.1&3.0E-16\\
dixmaana&3000&7&$<$0.1&1.0E+00&11&$<$0.1&1.0E+00&5&$<$0.1&1.0E+00\\
dixmaanb&3000&7&1.1E-01&1.0E+00&11&$<$0.1&1.0E+00&4&$<$0.1&1.0E+00\\
dixmaanc&3000&8&1.0E-01&1.0E+00&13&$<$0.1&1.0E+00&5&$<$0.1&1.0E+00\\
dixmaand&3000&10&1.7E-01&1.0E+00&16&$<$0.1&1.0E+00&5&$<$0.1&1.0E+00\\
dixmaane&3000&255&3.2E-01&1.0E+00&305&4.7E-01&1.0E+00&256&3.0E-01&1.0E+00\\
dixmaanf&3000&195&6.1E-01&1.0E+00&248&8.4E-01&1.0E+00&268&7.0E-01&1.0E+00\\
dixmaang&3000&227&7.0E-01&1.0E+00&239&8.5E-01&1.0E+00&268&7.2E-01&1.0E+00\\
dixmaanh&3000&181&5.6E-01&1.0E+00&316&1.0E+00&1.0E+00&220&6.0E-01&1.0E+00\\
 dixmaani&3000&6084&1.1E+01&1.0E+00&4793&7.2E+00&1.0E+00&4806&5.0E+00&1.0E+00\\
 dixmaanj&3000&3816&1.2E+01&1.0E+00&1409&4.6E+00&1.0E+00&675&1.7E+00&1.0E+00\\
 dixmaank&3000&3597&1.1E+01&1.0E+00&663&2.3E+00&1.0E+00&499&1.3E+00&1.0E+00\\
 dixmaanl&3000&2050&6.3E+00&1.0E+00&709&2.3E+00&1.0E+00&380&1.0E+00&1.0E+00\\
dixon3dq*&10&10&$<$0.1&2.2E-28&10&$<$0.1&2.2E-28&10&$<$0.1&2.4E-28\\
dqdrtic*&5000&6&$<$0.1&2.3E-15&6&$<$0.1&2.3E-15&6&$<$0.1&2.3E-15\\
dqrtic&5000&13&5.7E-01&1.0E-01&164&4.2E-01&9.5E-02&7&2.5E-01&9.5E-03\\
edensch&2000&17&$<$0.1&1.2E+04&17&$<$0.1&1.2E+04&16&$<$0.1&1.2E+04\\
eg2&1000&2&$<$0.1&-1.0E+03&2&$<$0.1&-1.0E+03&2&$<$0.1&-1.0E+03\\
engval1&5000&14&2.8E-01&5.5E+03&21&2.2E-01&5.5E+03&8&1.3E-01&5.5E+03\\
engval2&3&28&$<$0.1&5.9E-13&36&$<$0.1&9.0E-17&29&$<$0.1&1.5E-13\\
errinros&50&259&$<$0.1&4.0E+01&432&$<$0.1&4.0E+01&394&$<$0.1&4.0E+01\\
expfit&2&11&$<$0.1&2.4E-01&12&$<$0.1&2.4E-01&5&$<$0.1&2.4E-01\\
fletcbv2&100&97&$<$0.1&-5.1E-01&97&$<$0.1&-5.1E-01&97&$<$0.1&-5.1E-01\\
fletchcr&100&302&$<$0.1&3.5E-13&293&$<$0.1&1.3E-13&247&$<$0.1&2.4E-13\\
flosp2hl&650&(IL)&&&(IL)&&&(IL)&&\\
flosp2hm&650&(IL)&&&(IL)&&&(IL)&&\\
flosp2th&650&(IL)&&&(IL)&&&(IL)&&\\
flosp2tl&650&(IL)&&&(IL)&&&(IL)&&\\
flosp2tm&650&(IL)&&&(IL)&&&(IL)&&\\
fminsrf2&1024&223&1.9E-01&1.0E+00&384&3.2E-01&1.0E+00&1106&1.4E+00&1.0E+00\\
fminsurf&1024&202&1.6E-01&1.0E+00&339&3.5E-01&1.0E+00&420&5.8E-01&1.0E+00\\
freuroth&5000&22&7.3E-01&6.1E+05&28&3.1E-01&6.1E+05&54&7.9E+00&6.1E+05\\
genhumps&5&36&$<$0.1&7.0E-14&34&$<$0.1&5.5E-13&24&$<$0.1&3.5E-12\\
genrose&500&2668&6.8E-01&1.0E+00&3682&5.5E-01&1.0E+00&2572&9.5E-01&1.0E+00\\
growth&3&196&$<$0.1&1.0E+00&159&$<$0.1&1.0E+00&(IL)&&\\
growthls&3&169&$<$0.1&1.0E+00&177&$<$0.1&1.0E+00&(IL)&&\\
gulf&3&41&$<$0.1&4.7E-13&32&$<$0.1&4.7E-09&41&$<$0.1&4.9E-10\\
hairy&2&18&$<$0.1&2.0E+01&11&$<$0.1&2.0E+01&5&$<$0.1&2.0E+01\\
hatfldd&3&22&$<$0.1&2.5E-07&26&$<$0.1&2.5E-07&36&$<$0.1&2.6E-07\\
hatflde&3&35&$<$0.1&4.4E-07&28&$<$0.1&4.4E-07&44&$<$0.1&4.4E-07\\
heart6ls&6&(IL)&&&(IL)&&&(IL)&&\\
heart8ls&8&339&$<$0.1&1.5E-16&382&$<$0.1&1.5E-12&590&$<$0.1&5.2E-12\\
helix&3&21&$<$0.1&3.6E-17&25&$<$0.1&2.2E-17&18&$<$0.1&8.4E-23\\
\hline
  \end{tabular}}
    \caption{Table \ref{tab:cutestResults1} (Continued)}
    \label{tab:cutestResults2}
\end{table}

\begin{table}[ht]
\centering
\resizebox{\linewidth}{!}{
\begin{tabular}{|l|r|r|r|r|r|r|r|r|r|r| }
 \hline
  \multicolumn{2}{|c|}{} & \multicolumn{3}{|c|}{With Powell Restarts} &\multicolumn{3}{|c|}{No Powell Restarts} &\multicolumn{3}{|c|}{Hybrid Cubic} \\
 \hline
 Name&n&Iter&Time&f(x*)&Iter&Time&f(x*)&Iter&Time&f(x*)\\
 \hline
 hilberta*&10&8&$<$0.1&5.8E-10&8&$<$0.1&5.8E-10&8&$<$0.1&5.8E-10\\
 hilbertb*&50&5&$<$0.1&2.1E-20&5&$<$0.1&2.1E-20&5&$<$0.1&2.1E-20\\ 
himmelbb&2&4&$<$0.1&1.7E-18&5&$<$0.1&1.7E-16&4&$<$0.1&5.8E-19\\

 himmelbf&4&44&$<$0.1&3.2E+02&81&$<$0.1&3.2E+02&30&$<$0.1&3.2E+02\\
himmelbg&2&6&$<$0.1&1.6E-16&7&$<$0.1&9.4E-20&6&$<$0.1&1.3E-16\\
himmelbh&2&5&$<$0.1&-1.0E+00&5&$<$0.1&-1.0E+00&5&$<$0.1&-1.0E+00\\

 humps&2&55&$<$0.1&2.1E-12&89&$<$0.1&4.6E-16&48&$<$0.1&2.2E-12\\
jensmp&2&15&$<$0.1&1.2E+02&16&$<$0.1&1.2E+02&6&$<$0.1&1.2E+02\\
kowosb&4&25&$<$0.1&3.1E-04&45&$<$0.1&3.1E-04&63&$<$0.1&3.1E-04\\
liarwhd&10000&14&1.8E+00&5.5E-15&18&3.3E-01&1.3E-17&12&7.1E-01&2.4E-19\\
loghairy&2&184&$<$0.1&1.8E-01&66&$<$0.1&1.8E-01&4&$<$0.1&6.2E+00\\
mancino&100&15&1.6E-01&3.6E-15&15&1.9E-01&5.5E-15&15&1.8E-01&3.9E-15\\
maratosb&2&(E)&&&(E)&&&4&$<$0.1&-1.0E+00\\
methanb8&31&2101&$<$0.1&4.9E-05&1377&$<$0.1&5.2E-05&2602&1.5E-01&6.9E-05\\
methanl8&31&(IL)&&&(IL)&&&(IL)&&\\
mexhat&2&9&$<$0.1&-4.0E-02&6&$<$0.1&-4.0E-02&6&$<$0.1&-4.0E-02\\
meyer3&3&(E)&&&(E)&&&(IL)&&\\
minsurf&36&13&$<$0.1&1.0E+00&17&$<$0.1&1.0E+00&15&$<$0.1&1.0E+00\\
msqrtals&1024&3376&2.2E+01&2.7E-08&3658&2.9E+01&1.1E-08&3866&4.0E+01&4.7E-08\\
msqrtbls&1024&2376&1.6E+01&3.4E-09&2741&2.1E+01&1.7E-08&3630&3.8E+01&6.9E-09\\
nasty*&2&(E)&&&(E)&&&(E)&&\\
ncb20&1010&(E)&&&(E)&&&110&6.9E+00&1.7E+03\\
ncb20b&1000&(E)&&&(E)&&&34&7.4E+00&1.7E+03\\
nlmsurf&15129&2467&1.3E+02&3.9E+01&3857&9.9E+01&3.9E+01&4582&1.0E+02&3.9E+01\\
noncvxu2&1000&2684&1.1E+00&2.3E+03&2814&9.3E-01&2.3E+03&3325&1.1E+00&2.3E+03\\
noncvxun&1000&22&$<$0.1&2.3E+03&20&$<$0.1&2.3E+03&8&$<$0.1&2.3E+03\\
 nondia&9999&9&1.2E+00&2.8E-15&8&2.8E-01&3.9E-24&6&4.8E-01&3.2E-12\\
nondquar&10000&348&1.8E+01&6.4E-05&1523&8.8E+00&6.8E-05&683&1.4E+01&9.5E-05\\
nonmsqrt&9&670&$<$0.1&7.5E-01&282&$<$0.1&7.5E-01&925&$<$0.1&7.5E-01\\
osbornea&5&258&$<$0.1&5.5E-05&289&$<$0.1&5.5E-05&(IL)&&\\
osborneb&11&162&$<$0.1&4.0E-02&143&$<$0.1&4.0E-02&175&$<$0.1&4.0E-02\\
palmer1c*&8&(IL)&&&(E)&&&(IL)&&\\
 palmer1d*&7&(E)&&&(E)&&&8357&5.0E-01&6.5E-01\\
palmer1e&8&(IL)&&&(IL)&&&4&$<$0.1&0.0E+00\\
palmer2c*&8&(IL)&&&(IL)&&&(IL)&&\\
palmer2e&8&(IL)&&&(IL)&&&(IL)&&\\
palmer3c*&8&5844&1.0E-01&2.0E-02&(IL)&&&(IL)&&\\
palmer3e&8&(IL)&&&9272&1.9E-01&5.1E-05&(IL)&&\\
palmer4c*&8&3319&$<$0.1&5.0E-02&4021&$<$0.1&5.1E-02&(IL)&&\\
palmer4e&8&4453&$<$0.1&1.5E-04&6286&1.3E-01&1.5E-04&(IL)&&\\
palmer5c*&6&6&$<$0.1&2.1E+00&6&$<$0.1&2.1E+00&6&$<$0.1&2.1E+00\\
palmer5d*&4&9&$<$0.1&8.7E+01&9&$<$0.1&8.7E+01&8&$<$0.1&8.7E+01\\
palmer6c*&8&2336&$<$0.1&2.1E-02&6634&1.0E-01&1.9E-02&(IL)&&\\
palmer7c*&8&8231&1.4E-01&6.2E-01&(IL)&&&(IL)&&\\
palmer8c*&8&1863&$<$0.1&1.6E-01&4321&$<$0.1&1.7E-01&(IL)&&\\
penalty1&1000&25&$<$0.1&9.7E-03&49&$<$0.1&9.7E-03&25&$<$0.1&9.7E-03\\
penalty2&100&86&$<$0.1&9.7E+04&(E)&&&63&$<$0.1&9.7E+04\\

\hline
  \end{tabular}}
    \caption{Table \ref{tab:cutestResults1} (Continued)}
        \label{tab:cutestResults3}
\end{table}

\begin{table}[ht]
\centering
\resizebox{\linewidth}{!}{
\begin{tabular}{|l|r|r|r|r|r|r|r|r|r|r| }
 \hline
  \multicolumn{2}{|c|}{} & \multicolumn{3}{|c|}{With Powell Restarts} &\multicolumn{3}{|c|}{No Powell Restarts} &\multicolumn{3}{|c|}{Hybrid Cubic} \\
 \hline
 Name&n&Iter&Time&f(x*)&Iter&Time&f(x*)&Iter&Time&f(x*)\\
 \hline
 penalty3&100&(E)&&&(E)&&&(IL)&&\\
 pfit1&3&40&$<$0.1&2.9E-04&38&$<$0.1&2.9E-04&26&$<$0.1&2.9E-04\\
pfit1ls&3&40&$<$0.1&2.9E-04&38&$<$0.1&2.9E-04&26&$<$0.1&2.9E-04\\
pfit2&3&44&$<$0.1&1.2E-02&50&$<$0.1&1.2E-02&26&$<$0.1&1.2E-02\\
pfit2ls&3&44&$<$0.1&1.2E-02&50&$<$0.1&1.2E-02&26&$<$0.1&1.2E-02\\
pfit3&3&66&$<$0.1&8.2E-02&56&$<$0.1&8.2E-02&20&$<$0.1&8.2E-02\\
pfit3ls&3&66&$<$0.1&8.2E-02&56&$<$0.1&8.2E-02&20&$<$0.1&8.2E-02\\
pfit4&3&(E)&&&62&$<$0.1&2.6E-01&19&$<$0.1&2.6E-01\\
pfit4ls&3&(E)&&&62&$<$0.1&2.6E-01&19&$<$0.1&2.6E-01\\
powellsg&4&78&$<$0.1&2.0E-10&67&$<$0.1&2.2E-11&70&$<$0.1&5.6E-10\\
power*&1000&(IL)&&&(IL)&&&(IL)&&\\
quartc&10000&14&2.4E+00&1.6E-01&150&9.9E-01&5.0E-01&6&7.4E-01&7.3E-01\\
rosenbr&2&27&$<$0.1&9.4E-18&24&$<$0.1&7.1E-17&16&$<$0.1&1.1E-16\\
s201*&2&2&$<$0.1&4.7E-27&2&$<$0.1&4.7E-27&2&$<$0.1&4.7E-27\\
s202&2&9&$<$0.1&4.9E+01&8&$<$0.1&4.9E+01&7&$<$0.1&4.9E+01\\
s204&2&5&$<$0.1&1.8E-01&5&$<$0.1&1.8E-01&5&$<$0.1&1.8E-01\\
s205&2&9&$<$0.1&1.1E-16&10&$<$0.1&1.0E-17&8&$<$0.1&5.0E-13\\
s206&2&4&$<$0.1&1.9E-16&5&$<$0.1&2.1E-19&5&$<$0.1&1.9E-24\\
s207&2&7&$<$0.1&2.4E-13&8&$<$0.1&4.3E-14&8&$<$0.1&2.2E-13\\
s208&2&27&$<$0.1&9.4E-18&24&$<$0.1&7.1E-17&16&$<$0.1&1.1E-16\\
s209&2&98&$<$0.1&1.2E-18&86&$<$0.1&8.1E-22&39&$<$0.1&1.8E-19\\
s210&2&389&$<$0.1&5.3E-21&372&$<$0.1&3.6E-20&148&$<$0.1&1.7E-18\\
s211&2&14&$<$0.1&2.7E-16&17&$<$0.1&1.1E-17&16&$<$0.1&2.1E-20\\
s212&2&9&$<$0.1&6.9E-22&12&$<$0.1&1.1E-25&8&$<$0.1&2.2E-15\\
s213&2&10&$<$0.1&1.6E-12&14&$<$0.1&1.1E-09&12&$<$0.1&1.1E-09\\
s240*&3&2&$<$0.1&3.8E-15&2&$<$0.1&3.8E-15&2&$<$0.1&3.8E-15\\
s243&3&9&$<$0.1&8.0E-01&9&$<$0.1&8.0E-01&9&$<$0.1&8.0E-01\\
s245&3&11&$<$0.1&1.7E-15&11&$<$0.1&2.7E-17&30&$<$0.1&6.1E-14\\
s246&3&17&$<$0.1&1.8E-16&18&$<$0.1&5.5E-18&18&$<$0.1&4.1E-20\\
s256&4&78&$<$0.1&2.0E-10&67&$<$0.1&2.2E-11&70&$<$0.1&5.6E-10\\
s258&4&48&$<$0.1&6.2E-13&27&$<$0.1&1.8E-15&24&$<$0.1&3.1E-17\\
s260&4&48&$<$0.1&6.2E-13&27&$<$0.1&1.8E-15&24&$<$0.1&3.1E-17\\
s261&4&27&$<$0.1&1.2E-09&37&$<$0.1&6.4E-10&59&$<$0.1&1.1E-09\\
s266&5&12&$<$0.1&1.0E+00&13&$<$0.1&1.0E+00&10&$<$0.1&1.0E+00\\
s267&5&75&$<$0.1&2.6E-03&68&$<$0.1&7.7E-09&163&$<$0.1&5.5E-07\\
s271*&6&6&$<$0.1&0.0E+00&6&$<$0.1&0.0E+00&6&$<$0.1&0.0E+00\\
s272&6&34&$<$0.1&5.7E-03&61&$<$0.1&5.7E-03&69&$<$0.1&5.7E-03\\
s272a&6&67&$<$0.1&3.4E-02&68&$<$0.1&3.4E-02&(IL)&&\\
s273&6&11&$<$0.1&5.3E-18&14&$<$0.1&6.3E-15&6&$<$0.1&1.0E-14\\
s274*&2&2&$<$0.1&2.6E-24&2&$<$0.1&2.6E-24&2&$<$0.1&2.6E-24\\
s275*&4&3&$<$0.1&6.0E-12&3&$<$0.1&6.0E-12&3&$<$0.1&6.0E-12\\
s276*&6&3&$<$0.1&1.5E-12&3&$<$0.1&1.5E-12&3&$<$0.1&1.5E-12\\
 s281a*&10&11&$<$0.1&2.0E-15&11&$<$0.1&2.0E-15&11&$<$0.1&1.3E-16\\

\hline
  \end{tabular}}
    \caption{Table \ref{tab:cutestResults1} (Continued)}
    \label{tab:cutestResults4}
\end{table}

\begin{table}[ht]
\centering
\resizebox{\linewidth}{!}{
\begin{tabular}{|l|r|r|r|r|r|r|r|r|r|r| }
 \hline
  \multicolumn{2}{|c|}{} & \multicolumn{3}{|c|}{With Powell Restarts} &\multicolumn{3}{|c|}{No Powell Restarts} &\multicolumn{3}{|c|}{Hybrid Cubic} \\
 \hline
 Name&n&Iter&Time&f(x*)&Iter&Time&f(x*)&Iter&Time&f(x*)\\
 \hline
s282&10&212&$<$0.1&2.7E-15&220&$<$0.1&1.2E-16&292&$<$0.1&1.3E-15\\
s283&10&52&$<$0.1&1.5E-09&117&$<$0.1&7.3E-09&49&$<$0.1&2.4E-09\\
s286&20&24&$<$0.1&6.3E-16&27&$<$0.1&7.2E-14&22&$<$0.1&1.7E-17\\
s287&20&54&$<$0.1&2.4E-17&36&$<$0.1&9.3E-16&30&$<$0.1&9.4E-15\\
s288&20&70&$<$0.1&3.2E-10&80&$<$0.1&4.0E-10&58&$<$0.1&6.9E-10\\
 s289&30&4&$<$0.1&0.0E+00&4&$<$0.1&0.0E+00&3&$<$0.1&0.0E+00\\
  s290*&2&2&$<$0.1&1.1E-31&2&$<$0.1&1.1E-31&2&$<$0.1&3.1E-32\\
s291*&10&10&$<$0.1&2.8E-33&10&$<$0.1&2.8E-33&10&$<$0.1&5.5E-33\\
s292*&30&28&$<$0.1&7.1E-15&28&$<$0.1&7.1E-15&28&$<$0.1&7.1E-15\\
s293*&50&39&$<$0.1&3.0E-15&39&$<$0.1&3.0E-15&39&$<$0.1&3.0E-15\\
s294&6&50&$<$0.1&2.2E-15&54&$<$0.1&2.2E-17&48&$<$0.1&1.2E-20\\
s295&10&86&$<$0.1&1.2E-15&90&$<$0.1&3.8E-18&81&$<$0.1&2.0E-14\\
s296&16&115&$<$0.1&5.7E-14&134&$<$0.1&7.6E-15&117&$<$0.1&1.4E-14\\
s297&30&203&$<$0.1&1.6E-13&280&$<$0.1&1.3E-14&190&$<$0.1&1.1E-14\\
s298&50&288&$<$0.1&1.2E-14&413&$<$0.1&7.6E-15&317&$<$0.1&1.7E-14\\
s299&100&590&$<$0.1&5.3E-14&809&$<$0.1&1.1E-14&618&$<$0.1&3.7E-14\\
s300*&20&20&$<$0.1&-2.0E+01&20&$<$0.1&-2.0E+01&20&$<$0.1&-2.0E+01\\
s301*&50&50&$<$0.1&-5.0E+01&50&$<$0.1&-5.0E+01&50&$<$0.1&-5.0E+01\\
s302*&100&100&$<$0.1&-1.0E+02&100&$<$0.1&-1.0E+02&100&$<$0.1&-1.0E+02\\
s303&20&15&$<$0.1&6.7E-30&18&$<$0.1&3.3E-19&12&$<$0.1&8.1E-16\\
s304&50&11&$<$0.1&7.0E-14&19&$<$0.1&6.6E-25&9&$<$0.1&3.3E-24\\
s305&100&13&$<$0.1&4.2E-23&30&$<$0.1&1.6E-15&14&$<$0.1&3.1E-28\\
s308&2&7&$<$0.1&7.7E-01&8&$<$0.1&7.7E-01&7&$<$0.1&7.7E-01\\
s309&2&6&$<$0.1&2.9E-01&7&$<$0.1&2.9E-01&7&$<$0.1&2.9E-01\\
s311&2&6&$<$0.1&1.7E-14&7&$<$0.1&2.1E-23&5&$<$0.1&1.1E-19\\
s312&2&33&$<$0.1&5.9E+00&26&$<$0.1&5.9E+00&17&$<$0.1&5.9E+00\\
s314&2&4&$<$0.1&1.7E-01&5&$<$0.1&1.7E-01&5&$<$0.1&1.7E-01\\
s333&3&(E)&&&(E)&&&3&$<$0.1&0.0E+00\\
s334&3&17&$<$0.1&8.2E-03&16&$<$0.1&8.2E-03&15&$<$0.1&8.2E-03\\
s350&4&25&$<$0.1&3.1E-04&45&$<$0.1&3.1E-04&63&$<$0.1&3.1E-04\\
s351&4&65&$<$0.1&3.2E+02&46&$<$0.1&3.2E+02&54&$<$0.1&3.2E+02\\
s352*&4&4&$<$0.1&9.0E+02&4&$<$0.1&9.0E+02&4&$<$0.1&9.0E+02\\
s370&6&72&$<$0.1&2.3E-03&80&$<$0.1&2.3E-03&98&$<$0.1&2.3E-03\\
s371&9&771&$<$0.1&1.8E-06&1838&$<$0.1&4.0E-06&3361&1.3E-01&5.2E-06\\
s379&11&159&$<$0.1&4.0E-02&159&$<$0.1&4.0E-02&151&$<$0.1&4.0E-02\\
s386*&2&2&$<$0.1&4.7E-27&2&$<$0.1&4.7E-27&2&$<$0.1&4.7E-27\\
sbrybnd&5000&(IL)&&&(IL)&&&(IL)&&\\
schmvett&10000&(E)&&&(E)&&&(E)&&\\
scosine&10000&(IL)&&&(IL)&&&(IL)&&\\
scurly10&10000&(IL)&&&(IL)&&&(IL)&&\\
scurly20&10000&(IL)&&&(IL)&&&(IL)&&\\
scurly30&10000&(IL)&&&(IL)&&&(IL)&&\\
sineval&2&49&$<$0.1&2.8E-23&50&$<$0.1&4.3E-22&35&$<$0.1&5.2E-23\\
sinquad&10000&194&3.0E+01&4.1E-11&3259&4.0E+01&8.0E-09&1049&7.5E+01&2.6E-05\\sisser&2&7&$<$0.1&2.9E-10&12&$<$0.1&7.7E-11&4&$<$0.1&2.7E-10\\
snail&2&20&$<$0.1&1.7E-14&69&$<$0.1&3.1E-23&78&$<$0.1&2.5E-22\\

 \hline
  \end{tabular}}
    \caption{Table \ref{tab:cutestResults1} (Continued)}
    \label{tab:cutestResults5}
\end{table}

\begin{table}[ht]
\centering
\resizebox{\linewidth}{!}{
\begin{tabular}{|l|r|r|r|r|r|r|r|r|r|r| }
 \hline
  \multicolumn{2}{|c|}{} & \multicolumn{3}{|c|}{With Powell Restarts} &\multicolumn{3}{|c|}{No Powell Restarts} &\multicolumn{3}{|c|}{Hybrid Cubic} \\
 \hline
 Name&n&Iter&Time&f(x*)&Iter&Time&f(x*)&Iter&Time&f(x*)\\
 \hline
 srosenbr&10000&27&4.2E+00&6.7E-15&27&3.4E-01&3.3E-16&109&2.1E+00&5.0E-12\\
testquad*&1000&(IL)&&&(IL)&&&(IL)&&\\
tointgss&10000&6&1.0E+00&1.0E+01&4&2.3E-01&1.0E+01&3&3.7E-01&1.0E+01\\
tquartic&10000&9&1.4E+00&7.1E-15&10&2.8E-01&4.2E-14&9&7.6E-01&1.8E-12\\

tridia*&10000&(IL)&&&4128&1.1E+01&4.8E-15&5118&1.5E+01&4.9E-15\\
vardim&100&10&$<$0.1&4.6E-17&21&$<$0.1&7.5E-18&3&$<$0.1&7.3E-26\\
vibrbeam&8&(IL)&&&(IL)&&&(IL)&&\\
watson&31&2248&1.0E-01&1.1E-08&4969&2.6E-01&9.5E-09&8919&1.3E+00&9.7E-09\\
woods&10000&55&5.3E+00&1.8E-15&104&8.9E-01&6.5E-10&44&7.7E-01&6.2E-11\\
yfitu&3&67&$<$0.1&6.7E-13&63&$<$0.1&6.8E-13&78&$<$0.1&4.1E-06\\
zangwil2*&2&1&$<$0.1&-1.8E+01&1&$<$0.1&-1.8E+01&1&$<$0.1&-1.8E+01\\ \hline
 \hline
  \end{tabular}}
    \caption{Table \ref{tab:cutestResults1} (Continued)}
    \label{tab:cutestResults6}
\end{table}
 
\end{document}